\documentclass[12pt,leqno,amssymb,psfig,eepic,figfonts,epsf]{amsart}

\usepackage{amssymb, amsmath}

\renewcommand\hat{\widehat}

\def\emty{\emptyset}

\def\diam{\mbox{\rm diam}}

\def\dis{\displaystyle}
\newcommand\ben{\begin{enumerate}}
\newcommand\een{\end{enumerate}}
\newcommand\bit{\begin{itemize}}
\newcommand\eit{\end{itemize}}
\def\endp{\hspace*{\fill}$\rule{.55em}{.55em}$ \smallskip}

\def\BBB{{\mathcal B}}

\def\DDD{{\mathcal D}}
\def\EEE{{\mathcal E}}

\def\GGG{{\mathcal G}}
\def\HHH{{\mathcal H}}

\def\th{theorem }

\def\Th{Theorem }

\def\qc{quasiconformal }

\def\homeo{homeomorphism }

\def\nbhd{neighborhood }

\def\Proof{{\noindent\sc Proof.} }

\def\al{\alpha}
\def\be{\beta}

\def\Te{\Theta}
\def\g{\gamma}
\def\G{\Gamma}
\def\si{\sigma}

\def\vp{\varphi}
\def\ep{\varepsilon}
\def\la{\lambda}
\def\La{\Lambda}
\def\De{\Delta}
\def\de{\delta}
\def\om{\omega}

\def\hal{\hat{\al}}
\def\hbe{\hat{\be}}

\def\hw{\hat{w}}

\def\R{\mbox{$\mathbb R$}}

\def\Z{\mbox{$\mathbb Z$}}

\def\N{\mbox{$\mathbb N$}}
\def\SS{\mbox{$\mathbb S$}}
\def\HH{\mbox{$\mathbb H$}}

\def\PP{\mbox{$\mathbb P$}}

\def\EE{\mbox{$\mathbb E$}}
\def\tP{\mbox{$\widetilde{\PP}$}}

\newcommand{\fdefeq}{\overset{\text{def.}}{=}}

\newtheorem{newthm}{Theorem}
\newtheorem{theorem}{Theorem}[section]
\newtheorem{lemma}[theorem]{Lemma}
\newtheorem{proposition}[theorem]{Proposition}
\newtheorem{corollary}[theorem]{Corollary}
\newtheorem{defthm}[theorem]{Definition et \th}
\newcommand{\REFEQN}[1] { \begin{equation}\label{#1} }
\newcommand{\ENDEQN}{\end{equation}}
\newcommand{\REFTHM}[1] { \begin{theorem}\label{#1} }
\newcommand{\ENDTHM}{\end{theorem}}
\newcommand{\REFNTH}[1] { \begin{newthm}\label{#1} }
\newcommand{\ENDNTH}{\end{newthm}}
\newcommand{\REFPROP}[1]{\begin{proposition}\label{#1} }
\newcommand{\ENDPROP}{\end{proposition} }
\newcommand{\REFLEM}[1]{\begin{lemma}\label{#1} }
\newcommand{\ENDLEM}{\end{lemma} }
\newcommand{\REFCOR}[1]{\begin{corollary}\label{#1} }
\newcommand{\ENDCOR}{\end{corollary} }
\newcommand{\REFDEFTHM}[1] { \begin{defthm}\label{#1} }
\newcommand{\ENDDEFTHM}{\end{defthm}}


\title{Harmonic measures versus quasiconformal measures for hyperbolic groups}

\address{Eurandom, P.O. Box 513, 5600 MB Eindhoven, The Netherlands\\
LATP/CMI\\
Universit\'e de Provence\\
39 rue Fr\'ed\'eric Joliot-Curie\\
13453 Marseille Cedex 13, France}

\author{S\'ebastien Blach\`ere, Peter Ha\"{\i}ssinsky \& Pierre Mathieu}
\address{
LATP/CMI\\
Universit\'e de Provence\\
39 rue Fr\'ed\'eric Joliot-Curie\\
13453 Marseille Cedex 13, France}

\date{\today}

\begin{document}

\maketitle

\begin{abstract}
We establish a dimension formula for the harmonic measure of a finitely 
supported and symmetric random walk on a hyperbolic group. 
We also characterize random walks for which this dimension is maximal. 
Our approach is based on the Green metric, a metric which provides a geometric point of view on random walks
and, in particular, which allows us to interpret harmonic measures as \qc measures on the boundary
of the group.
\end{abstract}

\section{Introduction}

It is a leading thread in hyperbolic geometry to try to understand properties of hyperbolic spaces
by studying their large-scale behaviour. 
This principle is applied through the
introduction of a canonical compactification which characterises the space itself.
For instance a hyperbolic group $\G$ in the sense of Gromov admits
a natural boundary at infinity $\partial \G$: it is a topologically
well-defined
compact set on which $\G$ acts by homeomorphisms.
Together, the pair consisting of the boundary $\partial\G$ with the action of $\G$
characterises the hyperbolicity of the group. Topological properties of $\partial\G$ also encode the
algebraic structure of the group. For instance one proves that $\G$ is virtually free if and only if 
$\partial\G$ is a Cantor set (see \cite{sta}
and also 
\cite{bow} for other results in this vein).
Moreover, the boundary is endowed with a canonical \qc structure which
determines the quasi-isometry class of the group (see \cite{bk} and the references therein for details).

Characterising special subclasses of hyperbolic
groups such as cocompact Kleinian groups often requires the construction of
special metrics and measures on the boundary which carry some geometrical
information.
For example, M.\,Bonk and B.\,Kleiner proved that a group admits a cocompact Kleinian action
on the hyperbolic space $\HH^n$, $n\ge 3$, if and only if its boundary has topological dimension $n-1$ and carries
an Ahlfors-regular metric of dimension $n-1$ \cite{bkl}.

\bigskip

There are two main  constructions
of measures on the boundary of a hyperbolic group:
\qc measures and harmonic measures.
 Let us recall these constructions.

Given a cocompact properly discontinuous action of $\G$ by isometries on
a pointed proper geodesic metric space $(X,w,d)$, the Patterson-Sullivan procedure consists
in taking weak limits of
$$\frac{1}{\sum_{\g\in\G} e^{-s d(w,\g(w))}} \sum_{\g\in\G} e^{-sd(w,\g(w))}\de_{\g(w)}$$
as $s$ decreases to the  logarithmic
volume growth $$v\fdefeq\limsup_{R\to \infty}\frac{1}{R}\log |B(w,R)\cap \G( w)|\,.$$
Patterson-Sullivan measures are \qc measures and Hausdorff measures
of $\partial X$ when endowed with
a visual metric.

Given a probability measure $\mu$ on $\Gamma$, the random walk $(Z_n)_n$ starting from the neutral element
$e$ associated with $\mu$ is defined by
$$Z_0=e\,;\, Z_{n+1}=Z_n\cdot X_{n+1}\, , $$
where $(X_n)$ is a sequence of independent and identically distributed
random variables of law $\mu$.
Under some mild assumptions on $\mu$, the walk $(Z_n)_n$ almost surely converges to a point $Z_\infty\in\partial \G$.
The law
of $Z_\infty$ is by definition the harmonic measure $\nu$.

\bigskip

The purpose of this work is to investigate the interplay between those two classes of measures
and take advantage of this interplay to derive information on the geometry of harmonic measures.

The usual tool for this kind of results is to replace the action of the group by a linear-in-time action
of a dynamical system and then to apply the thermodynamic formalism to it: for free groups and Fuchsian groups,
a Markov-map $F_{\G}$ has been introduced on the boundary
which is orbit-equivalent to $\G$ \cite{bse,led}.
For discrete subgroups of isometries of a Cartan-Hadamard manifold,
one may work with the geodesic flow \cite{led2,led3,kai90, kai94}. 
Both these methods  seem difficult to implement for general hyperbolic groups.
On the one hand, it is not obvious how to associate a Markov map with a general hyperbolic group,
even using the automatic structure of the group.
On the other hand, the construction of the geodesic flow for general hyperbolic spaces is delicate
and its  mixing properties do not seem
strong enough to apply the thermodynamic formalism. 

In a different spirit, it is proved in \cite{connmuch} that any Patterson-Sullivan measure can be 
realized as the harmonic measure of some random walk. Unfortunately such a general statement without 
any information on the law of the increments of the random walk is not sufficient to provide any 
real insight in the behaviour of the walk.

Our approach directly combines geometric and probabilistic arguments.  
Since we avoid using the thermodynamic formalism, we believe it is more elementary. 
We make a heavy use of the so-called Green metric associated with the random walk, 
and we emphasize the connections between the geometry of this metric and
the properties of the random walk it comes from. 
The problem
of identifying the Martin and visual boundaries is an example of such a connection.
We give sufficient conditions  for the Green metric to be hyperbolic (although it is not geodesic in general).
On the other hand,  its explicit expression 
in terms of the hitting probability of the random walk makes it possible to directly 
take advantage of the independence of the increments of the walk. The combination of both 
facts yields very precise estimates on how random paths deviate from geodesics.

Thus we show  that, for a general hyperbolic group, 
the Hausdorff dimension of the harmonic measure
can be explicitly computed  and satisfies a 'dimension-entropy-rate of escape' formula 
and we characterise those harmonic measures of maximal dimension. 
Our point of view also allows us to get an alternative and rather straightforward 
proof of the fact that the 
harmonic measure of a random walk on a Fuchsian group with cusps is singular, a result 
previously established in \cite{guiv2} and \cite{deroin} by completely different methods.

The rest of this introduction is devoted to a more detailed description of our results.

\bigskip

\subsection{Geometric setting} \label{S1.1} Given a hyperbolic group $\G$, we let $\DDD(\G)$ denote
the collection of hyperbolic left-invariant metrics on $\G$ which are quasi-isometric to a word metric induced
by a finite generating set of $\G$.
In general these metrics do not come from proper geodesic metric spaces
as we will see (cf. \Th \ref{main1} for instance).
In the sequel, we will distinguish the group as a space and as acting on a space: we keep  the notation $\G$ for the group, and we
denote  by $X$  the group as a metric space endowed with a metric
$d\in\DDD(\G)$.
We may  equivalently write $(X,d)\in \DDD(\G)$. We will often require a base point which we will denote by $w\in X$.

This setting enables us to capture in particular the following two situations.
\bit
\item Assume that $\G$ admits a cocompact properly discontinuous action by isometries on a proper geodesic space $(Y,d)$. Pick $w\in Y$ such that $\g\in\G\mapsto \g(w)$ is a bijection, and
consider $X=\G(w)$ with the restriction of $d$.
\item We may choose $(X,d)=(\Gamma,d_G)$ where $d_G$ is the Green metric associated with a random walk (see \Th \ref{main1}).
\eit

Let $\mu$ be a symmetric probability measure the support of which generates $\G$.
Even if the support of $\mu$ may be infinite, we will require some compatibility
with the geometry of the quasi-isometry class of $\DDD(\G)$.
Thus, we will often assume one of the  following two  assumptions.
Given a metric $(X,d)\in\DDD(\G)$,
we say that the random walk has {\it finite first moment} if
$$\sum_{\g\in\G} d(w,\g(w))\mu(\g) <\infty\,.$$
We say that the random walk has an {\it exponential moment} if there exists $\la >0$ such that
$$\sum_{\g\in\G} e^{\la d(w,\g(w))}\mu(\g) <\infty\,.$$
Note that both these conditions only depend on the quasi-isometry class of the metric.

\subsection{The Green metric} \label{S1.2}
The analogy between both families of measures -- \qc and harmonic -- has already been pointed out in the literature.
Our first task is to make this empirical fact a theorem i.e., we prove that harmonic measures
are indeed \qc measures for a well-chosen metric: given a symmetric law $\mu$ on $\G$ such that
its support generates $\G$,
let $F(x,y)$ be the probability that the random walk started at $x$ ever hits
$y$. Up to a constant factor, $F(x,y)$ coincides with the Green function
$$
G(x,y)\fdefeq   \sum_{n=0}^{\infty} \PP^x [ Z_n=y ] =\sum_{n=0}^{\infty} \mu^n(x^{-1}y)\, ,$$
where $\PP^x$ denotes the probability law of the random walk $(Z_n)$ with $Z_0=x$
(if $Z_0=e$, the neutral element of $\G$, we will simply write $\PP^e=\PP)$,
and where, for each $n\ge 1$,  $\mu^n$ is the law of $Z_n$ i.e., the $n$th convolution power of the measure $\mu$.

We define the {\it Green  metric} between $x$ and $y$ in $\Gamma$ by
$$d_G(x,y)\fdefeq -\log F(x,y)\,.$$
This metric was
first introduced
by S.\,Blach\`ere and S.\,Brofferio in \cite{bb} and further studied in \cite{bhm}.
It is well-defined as soon as the walk is transient i.e., eventually leaves any finite set. This is the
case as soon as $\G$ is a non-elementary hyperbolic group.

Non-elementary hyperbolic groups are non-amenable and for such groups and finitely supported laws $\mu$, it was proved
in \cite{bb} that the Green and word  metrics are quasi-isometric.
Nevertheless it does not follow from this simple fact that $d_G$ is hyperbolic, see the discussion below,
\S\,\ref{intro:qrs}.

We first  prove the following:

\REFTHM{main1} Let $\Gamma$ be a non-elementary hyperbolic group,
$\mu$ a symmetric probability measure on $\Gamma$ the support of which
generates $\Gamma$.
\ben
\item[(i)] Assume that $\mu$ has an exponential moment, then $d_G\in\DDD(\G)$ if and only if
 for any $r$ there exists a constant $C(r)$ such that
\begin{equation}\label{eq:an} F(x,y)\le C(r)F(x,v)F(v,y)\end{equation}
whenever $x,y$ and $v$ are points in a locally finite Cayley graph of $\G$ and $v$ is at distance at most $r$
from a geodesic segment between $x$ and $y$.
\item[(ii)] If $d_G\in\DDD(\G)$ then the harmonic measure is Ahlfors regular of dimension $1/\ep$, when
$\partial \G$ is endowed with a visual metric $d_{\ep}^G$ of parameter $\ep> 0$ induced by $d_G$.
\een
\ENDTHM

Visual metrics are defined in the next section.

A.\,Ancona proved that (\ref{eq:an}) holds for finitely supported laws $\mu$.
Condition  (\ref{eq:an}) has also been
coined by V.\,Kaimanovich as the key ingredient in proving that the Martin boundary coincides with the geometric (hyperbolic)
boundary \cite[Thm 3.1]{kai94} (See also \S\,\ref{intro:martin} and \S\,\ref{sec:martin} for a further discussion on the relationships between
the Green metric and the Martin boundary).

\Th \ref{main1} in particular yields

\REFCOR{cmain1} Let $\Gamma$ be a non-elementary hyperbolic group,
$\mu$ a finitely supported symmetric probability measure on $\Gamma$ the support of which
generates $\Gamma$. Then its associated Green metric $d_G$ is a left-invariant hyperbolic
metric on $\Gamma$ quasi-isometric to $\G$ such that the harmonic measure is Ahlfors regular of dimension $1/\ep$, when
$\partial \G$ is endowed with a visual metric $d_{\ep}^G$ of parameter $\ep> 0$ induced by $d_G$.
\ENDCOR

Our second source of examples of random walks  satisfying (\ref{eq:an}) will come from Brownian motions on Riemannian manifolds of negative curvature.
The corresponding law $\mu$ will then
have infinite support (see \S\,\ref{intro:bm} and \S\,\ref{sec:bm}).

\subsection{Dimension of the harmonic measure at infinity} Let $(X,d)\in\DDD(\G)$.
We fix a base point $w\in X$ and consider the random walk on $X$ started at $w$ i.e., the
sequence of $X$-valued random variables $(Z_n(w))$ defined by the action of $\Gamma$ on $X$.
There are (at least) two natural asymptotic quantities one can consider: the {\it asymptotic
entropy} $$h \fdefeq \lim_n \frac{- \sum_{\gamma \in \Gamma} \mu^n (\gamma) \log \mu^n(\gamma)}{n} = \lim_n \frac{- \sum_{x \in \Gamma (w)} \PP [Z_n(w)=x]  \log \PP [Z_n(w)=x]}{n}
$$ which measures the way the law of $Z_n(w)$ is spread in different directions,
and the {\it rate of escape} or {\it drift}
$$\ell\fdefeq \lim_n \frac{d(w,Z_n(w))}{n}\, ,
$$
which estimates how far $Z_n(w)$ is from
its initial point $w$.  (The above limits for $h$ and $\ell$ are almost sure and in $L^1$ and they are finite as soon as the law has a finite  first moment.)

We obtain the following.

\REFTHM{main2} Let  $\Gamma$ be a non-elementary hyperbolic group, $(X,d)\in\DDD(\G)$,
$d_{\ep}$ be a visual metric of $\partial X$, and let
$B_\ep(a,r)$ be the ball of center $a\in\partial X$ and radius $r$ for the distance $d_\ep$.
Let $\nu$  be the harmonic measure of a random walk $(Z_n)$
whose increments are given by a symmetric law $\mu$ with finite  first moment such that $d_G\in\DDD(\G)$.

The pointwise Hausdorff dimension $\lim_{r\to 0} \frac{\log \nu (B_\ep(a,r))}{\log r}$ exists for $\nu$-almost every $a\in\partial X$,
and is independent from  the choice of $a$. More precisely, for $\nu$-almost every $a\in\partial X$,
$$ \lim_{r\to 0} \frac{\log \nu (B_\ep(a,r))}{\log r}= \frac{\ell_G}{\ep\ell}$$
where $\ell>0$ denotes the rate of escape of the walk with respect to $d$ and
$\ell_G\fdefeq \lim_n \frac{d_G(w,Z_n(w))}{n}$ the rate of escape with respect to $d_G$.
\ENDTHM

We recall that the dimension of a measure is the infimum Hausdorff dimension of sets of
positive measure.
In \cite{bhm}, it was shown that $\ell_G= h$ the asymptotic entropy of the walk.
 From \Th\ref{main2},
we deduce that

\REFCOR{cmain2} Under the assumptions of \Th \ref{main2},
  $$\dim\nu = \frac{h}{\ep\ell}$$
where $h$ denotes the asymptotic entropy of the walk and $\ell$ its
rate of escape with respect to $d$.\ENDCOR

This dimension formula already appears in the work of F.\,Ledrappier for random walks on free groups \cite{led}.
See also V.\,Kaimanovich, \cite{ka2}.
For general hyperbolic groups, V.\,Leprince established the inequality $\dim\nu\le h/(\ep\ell)$
and made constructions of harmonic measures with arbitrarily small dimension \cite{lep}.
More recently,  V.\,Leprince established that $h/\ep\ell$ is also the box dimension of the harmonic
measure under the sole assumption that the random walk has a finite first moment \cite{lep2}.
Note however that the notion of box dimension
is too weak to ensure the existence of the pointwise Hausdorff dimension almost everywhere.

This formula is also closely related to the dimension formula proved
for ergodic invariant measures with positive entropy in the context of geometric dynamical systems:
the drift corresponds to a Lyapunov exponent \cite{lsy}.

\bigskip

\subsection{Characterisation of harmonic measures with maximal dimension}
Given a random walk  on a finitely generated group $\G$ endowed with a left-invariant
metric $d$, the so-called {\it fundamental inequality} between the asymptotic entropy $h$, the drift $\ell$
and  the logarithmic growth rate $v$ of the action of
$\G$ reads 
$$h\le \ell v\,.$$  It holds as soon as all these objects are well-defined (cf. \cite{bhm}).
Corollary \ref{cmain2}
provides a geometric interpretation of this inequality in terms of the harmonic measure:
indeed, since $v/\ep$ is the dimension of $(\partial X,d_{\ep})$, see \cite{co}, it is clearly
larger than the dimension of $\nu$. 

A.\,Vershik suggested the study of the case of equality (see \cite{yg,av}). For any hyperbolic group,
\Th \ref{main1} implies that the equality $h= \ell v$ holds for the Green metric and \Th \ref{main3}
below 
shows
that the equality for some $d\in\DDD(\G)$
implies  $d$ is almost proportional to $d_G$. In particular,
given a metric in $\DDD(\G)$,  all the  harmonic measures for which the (fundamental) equality
holds belong to the same class of \qc measures.

\bigskip
In the sequel, two measures will be called {\it equivalent} if they share the same sets of zero measure.

\REFTHM{main3}  Let  $\Gamma$ be a non-elementary hyperbolic group and
$(X,d)\in\DDD(\G)$; let $d_{\ep}$ be a visual metric of $\partial X$, and $\nu$ the harmonic measure given by a symmetric
law $\mu$ with an exponential moment, the support of
 which generates $\Gamma$.
We further assume that $(X,d_G)\in\DDD(\G)$. 
We denote by $\rho$ a \qc measure on
$(\partial X,d_\ep)$. The following propositions are equivalent.
\bit
\item[(i)] We have the equality $h=\ell v$.
\item[(ii)] The measures $\rho$ and $\nu$ are equivalent.
\item[(iii)] The measures $\rho$ and $\nu$ are equivalent
and the density is almost surely bounded and
bounded away from $0$.
\item[(iv)] The map $(\G, d_G)\stackrel{Id}{\longrightarrow} (X,vd)$ is a $(1,C)$-quasi-isometry.
\item[(v)] The measure $\nu$ is a \qc measure of $(\partial X,d_\ep)$ .
\eit\ENDTHM

This \th is the counterpart of a result of F.\,Ledrappier for Brownian motions on universal 
covers of compact Riemannian manifolds
of negative sectional curvature \cite{led2}, see also \S\,\ref{intro:bm}.
Similar results have been established for the free group
 with free generators, see \cite{led}.
The case of equality $h=\ell v$ has also been studied for particular sets of generators of
free products of finite groups \cite{maimat}. For universal covers of finite graphs, see \cite{russ}.

\Th \ref{main3} enables us to compare random walks and decide when their harmonic measures are equivalent.

\REFCOR{cmain3}
Let $\Gamma$ be a non-elementary hyperbolic group with two finitely supported symmetric  probability measures $\mu$ and $\hat{\mu}$ where both supports generate $\Gamma$.
We consider the random walks $(Z_n)$ and $(\hat{Z}_n)$.
Let us denote their Green functions by $G$ and $\hat{G}$ respectively,
the asymptotic entropies by $h$ and $\hat{h}$,
and the harmonic measures seen from the neutral element $e$ by $\nu$ and $\hat{\nu}$.
The following propositions are equivalent.
\bit
\item[(i)] We have the equality $$\hat{h}=\lim \frac{-1}{n} \log G(e,\hat{Z_n}) $$ in $L^1$ and almost surely.
\item[(ii)] We have the equality $$h=\lim \frac{-1}{n} \log \hat{G}(e,Z_n) $$ in $L^1$ and almost surely.
\item[(iii)] The measures $\nu$ and $\hat{\nu}$ are equivalent.
\item[(iv)] There is a constant $C$ such that
$$\frac{1}{C} \le \frac{G(x,y)}{\hat{G}(x,y)}\le  C\,.$$
\eit
\ENDCOR

\bigskip

\subsection{The Green metric and the Martin compactification} \label{intro:martin}
Given a probability measure $\mu$ on a countable group $\G$,  one defines the Martin kernel
$$
K(x,y) =K_y(x)\fdefeq \frac{G(x,y)}{G(e,y)}\,.
$$
By definition, the {\it Martin compactification} $\Gamma \cup \partial_M \Gamma$
is the smallest compactification
of $\G$ endowed with the discrete topology such that the  Martin kernel
continuously extends to $\G\times(\G\cup\partial_M\G)$. Then
 $\partial_M \Gamma$ is called the {\it Martin boundary}.

A general theme is to identify the Martin boundary with a geometric boundary of the group. It was observed in \cite{bhm} that
the Martin compactification coincides with the Busemann compactification of $(\G,d_G)$. We go one step further by showing that the
Green metric provides a common framework for the identification of the Martin boundary with the boundary at infinity of a hyperbolic space
(cf. \cite{aa,bl,kai94}).

\REFTHM{martin} Let $\G$ be a countable group, $\mu$ a symmetric
probability measure the support of which generates $\G$.
We assume that the corresponding random walk is transient.
If the Green metric is hyperbolic, then the Martin boundary
consists only of minimal points and it is homeomorphic to
the hyperbolic boundary of $(\G,d_G)$.

In particular, if $\G$ is a non-elementary hyperbolic group
and if $d_G\in\DDD(\G)$, then $\partial_M\G$ is homeomorphic
to $\partial\G$.\ENDTHM

One easily deduces  from Corollary \ref{cmain1}:

\REFCOR{cor1:martin} {\em\bf (A.\,Ancona)}  Let $\G$ be a non-elementary
hyperbolic group, $\mu$ a finitely supported probability
measure the support of which generates $\G$. Then the Martin
boundary is homeomorphic to the hyperbolic boundary of $\G$.
\ENDCOR

In \S\,\ref{exples}, we provide examples of hyperbolic groups with random walks for which
the Green metric is hyperbolic, but not in the quasi-isometry class of the group, and also examples
of non-hyperbolic groups for which the Green metric is nonetheless hyperbolic. These examples
are constructed by discretising Brownian motions on Riemannian manifolds (see below).

\bigskip

\subsection{Brownian motion revisited} \label{intro:bm}
Let $M$ be the universal covering of a compact Riemannian manifold of negative curvature
with deck transformation
group $\G$ i.e., the action of $\G$ is isometric, cocompact and  properly discontinuous. The Brownian motion $(\xi_t)$ on $M$
is the diffusion process generated by the Laplace-Beltrami operator. It is known that the Brownian motion
 trajectory almost surely converges to some limit point  $\xi_\infty\in\partial M$. The law of $\xi_\infty$ is the harmonic measure of
the Brownian motion. The notions of rate of escape and asymptotic
entropy also make perfect sense in this setting.

Refining a method of T.\,Lyons and D.\,Sullivan \cite{ls}, W.\,Ballmann and F.\,Ledrappier construct in \cite{bl} a random
walk on $\G$ which mirrors the trajectories of the Brownian motion and to which we may apply our previous results.
This enables us to recover the following results.

\REFTHM{main4} Let $M$ be the universal covering of a compact Riemannian manifold of negative curvature
with logarithmic volume growth $v$. Let $d_\ep$ be a visual distance on $\partial M$. Then
$$\dim \nu = \frac{h_M}{\ep\ell_M}$$ where $h_M$ and
$\ell_M$ denote the asymptotic entropy and the drift of the Brownian motion respectively. Furthermore, $h_M=\ell_M v$
if and only if $\nu$ is equivalent to the Hausdorff measure of dimension $v/\ep$ on $(\partial M,d_\ep)$.\ENDTHM

The first result is folklore and explicitely stated by V.\,Kaimanovich in the introduction of \cite{kai90}, but we know
of no published proof. The second statement is due to F.\,Ledrappier \cite{led2}. Note that more is known: the equality
$h_M=\ell_M v$ is equivalent to the equality of $\nu$ with the canonical {\it conformal} measure on  $(\partial M,d_\ep)$,
and this is possible only if  $M$ is a rank $1$ symmetric space \cite{led3,bcg}.
\bigskip

\subsection{Quasiruled hyperbolic spaces}\label{intro:qrs} As previously mentioned,  S.\,Blach\`ere and
S.\,Brofferio proved that,  for finitely supported laws, the Green metric $d_G$ is quasi-isometric to the word metric.
But since $d_G$ is defined only on a countable set, it is unlikely to be the restriction of a proper geodesic
metric (which would have guaranteed the hyperbolicity of $(\G,d_G)$).
Therefore, the proof of \Th \ref{main1} requires the understanding of which metric spaces among the
quasi-isometry class of a given geodesic hyperbolic space are also hyperbolic.
For this, we coin the
notion of a quasiruler: a {\it $\tau$-quasiruler} is a quasigeodesic
$g:\R\to X$ such that, for
any $s<t<u$, $$d(g(s),g(t))+d(g(t),g(u))- d(g(s),g(u)) \le 2\tau.$$
 A metric space will be {\it quasiruled} if constants $(\la,c,\tau)$ exist so that the space is
$(\la,c)$ - quasigeodesic and if every $(\la,c)$-quasigeodesic
is
a $\tau$-quasiruler. 
We refer to the Appendix for details on the definitions and properties of quasigeodesics and quasiruled spaces. 
We prove the following characterisation of hyperbolicity, interesting in its own right.

\REFTHM{hypqrule} Let $X$ be a geodesic hyperbolic metric space, and $\vp:X\to Y$ a quasi-isometry, where $Y$ is
a metric space. Then $Y$ is hyperbolic if and only if it is quasiruled.
\ENDTHM

Theorem \ref{hypqrule} will be used to prove  that the hyperbolicity of $d_G$ 
is equivalent to condition (\ref{eq:an}) in Theorem \ref{main1}. 
We complete this discussion by exhibiting
for any hyperbolic group,  a non-hyperbolic left-invariant metric in its quasi-isometry class 
(cf. Proposition \ref{nonhyp}).

\subsection{Fuchsian groups with cusps} 
We provide an alternative proof based on the Green metric of a  theorem due to Y.\,Guivarc'h and Y.\,Le Jan 
about random walks on Fuchsian groups, see the last corollary of \cite{guiv2}.

\REFTHM{GuivLeJan} {\em\bf (Y.\,Guivarc'h \& Y.\,Le Jan)} 
Let $\Gamma$ be a discrete subgroup of $PSL_2(\R)$  
such that the quotient space $\HH^2/\Gamma$ is not compact but has finite volume. 
Let $\nu_1$ be the harmonic measure on $\SS^1$ given by a symmetric law $\mu$ with finite support. 
(Almost any trajectory of the random walk converges to a 
point in $\partial\HH^2= \SS^1$ and $\nu_1$ is the law of this limit point.) 
Then $\nu_1$ is singular with respect to the Lebesgue measure on $\SS^1$. 
\ENDTHM

Note that it follows from Theorem \ref{main3} that $\nu_1$ is singular with respect to the Lebesgue measure 
if and only if its dimension is less than $1$. 

This theorem was originally derived from results 
on winding numbers of the geodesic flow, see 
\cite{guiv2} and \cite{guiv}. 
A more recent proof based on ergodic properties of smooth group actions on $\SS^1$ 
was obtained by B.\,Deroin, V.\,Kleptsyn and A.\,Navas in \cite{deroin}. 
It applies to random walks with a finite first moment. 

We shall see how Theorem \ref{GuivLeJan} can also be deduced from the hyperbolicity 
of the Green metric through a rather straighforward argument. 
We only consider the symmetric and finite support case 
even though it would also work if the random 
walk has a first finite moment and if $d_G\in\DDD(\Gamma)$.

We thank B.\,Deroin, Y.\,Guivarc'h and Y.\,Le Jan for enlightening explanations on their theorem.

\subsection{Outline of the paper} In Section 2, we recall the main facts on hyperbolic groups which will
be used in the paper. In Section 3, we recall the construction of random walks,
discuss some of their properties and introduce the Green metric. We also prove  \Th \ref{martin} and \Th \ref{main1}.
We then draw some consequences on the harmonic
measure and the random walk.  The following Section 4 deals with the
proof of \Th \ref{main2}.  In Section 5, we deal with \Th \ref{main3} and its corollary and we conclude 
with the proof of Theorem \ref{GuivLeJan}.
Finally, Theorem \ref{main4} is proved in Section 6.
The appendices are devoted to quasiruled spaces. We prove \Th \ref{hypqrule} in Appendix A, and we show that
quasiruled spaces retain most properties of geodesic hyperbolic spaces: in Appendix B, we show
that the approximation of finite configurations
by trees still hold, and  we explain why M.\,Coornaert's \th on quasiconformal
measures remains valid in this setting.

\subsection{Notation} A distance in a metric space will be denoted
either by $d(\cdot,\cdot)$ or $|\cdot-\cdot|$. If $a$ and $b$ are positive,
$a\lesssim b$  means that there is a universal positive constant
$u$ such that $a\le ub$. We will write $a\asymp b$ when both $a\lesssim b$ and
$b\lesssim a$ hold.
Throughout the article, dependance of a constant on structural parameters of the space will not be notified unless needed.
Sometimes, it will be convenient to use Landau's notation $O(\cdot)$.


\section{Hyperbolicity in metric spaces}\label{S2}

Let $(X,d)$ be a metric space.
It is said to be {\it proper} if closed balls of finite radius are compact.
A {\it geodesic curve (resp. ray, segment)} is a curve isometric to $\R$
(resp. $\R_+$, a compact interval of $\R$).
The space $X$ is said to be {\it geodesic} if every pair of points can be
joined by a geodesic segment.

\bigskip

Given three points $x,y,w\in X$, one defines the Gromov inner product
as follows:

$$(x|y)_w \fdefeq (1/2)\{|x-w|+|y-w|-|x-y|\}\,.$$

\noindent{\bf Definition.} A metric space $(X,d)$ is  {\it $\de$-hyperbolic}
($\de\ge 0$) if, for any $w,x,y,z\in X$, the following ultrametric type inequality holds
$$(y|z)_w \ge \min \{(x|y)_w, (x|z)_w \} -\de\,.$$

  We shall write
$(\cdot|\cdot)_w=(\cdot|\cdot)$ when the choice of $w$ is clear from the context.

Hyperbolicity is a large-scale property of the space. To capture
this information, one defines the notion of quasi-isometry.

\noindent{\bf Definition.}  Let $X,Y$ be two metric spaces and $\la\ge 1$, $c\ge 0$ two constants.
A map
$f:X\to Y$ is a {\it $(\la,c)$-quasi-isometric embedding} if, for any
$x,x'\in X$, we have $$\frac{1}{\la} |x-x'| -c\le|f(x)-f(x')|\le \la|x-x'|+c\,.$$
The map $f$ is a {\it $(\la,c)$-quasi-isometry} if,  in addition, there exist a quasi-isometric embedding
$g:Y\to X$ and a constant $C$ such that  $|g\circ f(x)-x|\le C$ for any $x\in X$.
Equivalently, $f$ is a quasi-isometry if it is a quasi-isometric embedding such that
$Y$ is contained in a $C$-\nbhd of $f(X)$. We then say that $f$ is {\it $C$-cobounded}.

In the sequel, we will always choose the constants so that that a $(\la,c)$-quasi-isometry is $c$-cobounded.

\bigskip

\noindent{\bf Definition.} A {\it quasigeodesic curve (resp. ray, segment)} is the image of
$\R$ (resp. $\R_+$, a compact interval of $\R$) by a quasi-isometric embedding.

\bigskip

In a geodesic hyperbolic metric space $(X,d)$, quasigeodesics always shadow genuine geodesics i.e.,  given
a $(\la,c)$-quasigeodesic $q$, there is  a geodesic $g$ such that $d_H(g,q)\le K$, where
$d_H$ denotes the Hausdorff distance, and $K$ only depends on $\de$, $\la$ and $c$ \cite[Th. 5.6]{gh}.

\bigskip

{\noindent\bf Compactification.} Let $X$ be a proper hyperbolic space, and $w\in X$ a base point.
A sequence $(x_n)$ tends to infinity if, by definition, $(x_n|x_m)\to \infty$ as $m,n\to \infty$.
The  {\it visual} or {\it hyperbolic} boundary $\partial X$ of $X$ is the set of sequences which tend to infinity modulo the equivalence
relation defined by:  $(x_n)\sim (y_n)$ if $(x_n|y_n)\to\infty$.
One may also extend the Gromov inner product to points at infinity  in such a way that the inequality
$$(y|z) \ge \min \{(x|y), (x|z) \} -\de\,,$$ now holds for any
points $w,x,y,z\in X\cup\partial X$.

For each $\ep>0$ small enough, there exists a so-called {\it visual metric} $d_{\ep}$ on $\partial X$ i.e  which satisfies for any $a,b\in\partial X$:
 $d_{\ep}(a,b)\asymp e^{-\ep (a|b)}$.

We shall use the notation $B_{\ep}(a,r)$ to denote the ball in the space $(\partial X,d_{\ep})$ with center $a$
and radius $r$.

We refer to \cite{gh} for the details (chap. 6 and 7).

\bigskip

{\noindent\bf Busemann functions.} Let us assume that $(X,d)$ is a hyperbolic space.
Let $a\in\partial X$, $x,y\in X$.
 The function
$$\be_a(x,y) \fdefeq \sup\left\{\limsup_{n\to\infty} [d(x,a_n)-d(y,a_n)]\right\}\,,$$
where the supremum is taken over all sequences $(a_n)_n$ in $X$ which tends to $a$,
is called the {\it Busemann function} at the point $a$.

\bigskip

{\noindent\bf Shadows.} Let $R>0$ and $x\in X$. The    {\it  shadow} $\mho(x,R)$ is the set of points $a\in\partial X$ such
that $(a|x)_w\ge d(w,x)-R$.

Approximating finitely many points by a tree (cf. \Th \ref{aptrees}) yields:
\REFPROP{shadowVSball}
Let $(X,d)$ be a   hyperbolic space.
For any $\tau \ge 0$, there exist positive constants $C,R_0$ such that for any $R>R_0$, $a\in \partial X$ and $x\in X$ such
that $(w|a)_x\le \tau$,

$$
B_\ep\left( a, \frac{1}{C} e^{R\ep}e^{-\ep |w-x|} \right) \subset \mho (x,R) \subset B_\ep \left( a, C e^{R\ep}e^{-\ep |w-x|} \right)  \, .
$$
\ENDPROP

Shadows will enable us to control measures
on the boundary of a hyperbolic group, see the
lemma of the shadow in the next paragraph.

\subsection{Hyperbolic groups}
Let $X$ be a hyperbolic proper metric space and $\G$ a subgroup
of isometries which acts properly discontinuously on $X$ i.e.,  for any compact sets $K$ and $L$, the number of
group elements $\g\in\G$ such that $\g(K)\cap L\ne\emty$ is finite.
For any point $x\in X$, its orbit $\G(x)$ accumulates only  on the boundary $\partial X$, and its
set of accumulation points turns out to be independent of the choice of $x$;
by definition, $\overline{\G(x)}\cap\partial X$ is the {\it limit set} $\La(\G)$ of $\G$.

An action of a group $\Gamma$ on a metric space is said to be {\it geometric} if
\ben
\item each element acts by isometry;
\item the action is properly discontinuous;
\item the action is cocompact.
\een

For example, if $\G$ is a finitely generated group, $S$ is a finite symmetric set of generators,
one may consider the Cayley graph $X$ associated with $S$: the set of vertices are the elements of
the group, and pairs $(\g,\g')\in\G\times\G$ define an edge if $\g^{-1}\g'\in S$. Endowing $X$ with the
metric which makes each edge isometric to the segment $[0,1]$ defines the {\it word metric associated with $S$}.
It turns $X$ into a geodesic proper metric space on which $\G$ acts geometrically by left-translation.

We recall \v{S}varc-Milnor's lemma which provides a sort of converse statement, see \cite{gh}:

\REFLEM{svmi} Let $X$ be a geodesic proper metric space, and $\G$ a group which acts geometrically
on $X$. Then $\G$ is finitely generated and $X$ is quasi-isometric to any locally finite Cayley graph
of $\G$.\ENDLEM

\bigskip

A group $\Gamma$ is {\it hyperbolic} if it acts geometrically on a geodesic proper
hyperbolic metric space (e.g. a locally finite Cayley graph).  In this case, one has $\La(\G)=\partial X$.
Then \v{S}varc-Milnor's lemma above implies that $\G$ is finitely generated.

We will say that a metric space
$(X,d)$ is {\it quasi-isometric} to the group $\G$ if
it is quasi-isometric to a locally finite Cayley graph of $\G$.

Let $\G$ be a hyperbolic group geometrically acting
on $(X,d)$.
The action of $\Gamma$ extends to the boundary.
Busemann functions, visual metrics and the action of $\G$ are related by the following property:  
for any $a\in\partial X$ and any $\g\in\G$, there
exists 
a \nbhd $V$ of $a$ such that, for any $b,c\in V$,
$$d_\ep(\g(b),\g(c))\asymp L_\g(a) d_\ep(b,c)\,$$  where $L_\g(a) \fdefeq e^{\ep\be_a(w,\g^{-1}(w))}$.
Moreover, $\Gamma$ also acts on measures on  $\partial X$
through the rule $\g^*\rho(A) \fdefeq \rho(\g A)$.

\bigskip

A hyperbolic group is said to be {\it elementary} if it is finite or quasi-isometric to $\Z$.
We will only be dealing with non-elementary hyperbolic groups.

\subsection{Quasiconformal measures}

We now assume that $\Gamma$ is a hyperbolic group acting on a proper 
quasiruled hyperbolic metric space $(X,d)$.

The next theorem summarizes the main properties of \qc measures on the boundary of $X$.
It was proved by  M.\,Coornaert in \cite{co} in the context of geodesic spaces. We state here a more general
version to cover the case $d\in\DDD(\G)$. We justify the validity of this generalisation
at the end of the appendix.
We refer to Section 4 for the definitions of the Hausdorff measure and dimension.

\bigskip

\REFTHM{reg2} \label{reg} 
Let $(X,d)$ be a proper quasiruled hyperbolic space endowed with a geometric action of a 
non-elementary hyperbolic group $\Gamma$.  
For any small enough $\ep>0$, we have $0<\mbox{\em dim}_H\,(\partial X,d_\ep)<\infty$ and
$$v \fdefeq \limsup \frac{1}{R}\log\left|\{\G(w)\cap B(w,R)\}\right|=\ep\cdot\mbox{\em dim}_H\,(\partial X,d_\ep)\,.$$
Let $\rho$ be the  Hausdorff measure on $\partial X$ of dimension $\al \fdefeq v/\ep$\,;
\bit
\item[(i)] $\rho$ is Ahlfors-regular of dimension $\al$ i.e., for any $a\in\partial X$,
for any $r\in (0,\diam\partial X)$, $\rho(B_{\ep}(a,r))\asymp r^\al$.
In particular, $0<\rho(\partial X)<\infty$.\\
\item[(ii)] $\rho$ is a quasiconformal measure i.e.,
for any isometry $\g$  we have $\rho\ll \g^*\rho\ll\rho$ and
$$\frac{ d \g^*\rho}{d\rho}\asymp (L_\g)^\al \, \rho\ a.e.\,.$$\\
\item[(iii)] The action of $\Gamma$ is ergodic for $\rho$ i.e.,
for any $\Gamma$-invariant Borelian $B$ of $\partial X$,
$$
\rho (B) =0 \;\; \mathrm{or} \;\; \rho (\partial X \backslash B)=0 \, .
$$\eit

Moreover, if $\rho'$ is another  $\Gamma$-quasiconformal measure, then $\rho\ll\rho'\ll\rho$ and
$\dis\frac{ d\rho}{d\rho'}\asymp 1$ a.e. and
$$\left|\{\G(w)\cap B(w,R)\}\right|\asymp e^{vR}\,.$$ \ENDTHM

The class of measures thus defined on $\partial X$ is called the
{\it Patterson-Sullivan} class.
It does not depend on the choice of the parameter $\ep$ but it does depend on the metric $d$.

The study of \qc measures yields the following key estimate \cite{co}:

\REFLEM{shadow} \label{ashadow} {\em\bf (Lemma of the shadow)} Under the assumptions of \Th\ref{reg2}, there  exists $R_0$, such that
if $R>R_0$, then, for any $x\in X$,  $$\rho(\mho(x,R))\asymp  e^{-v d(w,x)}$$
where the implicit constants do not depend on $x$. \ENDLEM


\section{Random walks and Green metrics for hyperbolic groups}

Let $\G$ be a hyperbolic group, and let us consider the set $\DDD(\G)$ of left-invariant hyperbolic metrics on $\G$ which
are quasi-isometric to $\G$.
We fix such a metric $(X,d)\in\DDD(\G)$ with a base point $w\in X$, and we consider
a symmetric probability measure $\mu$
on $\G$ with finite  first moment i.e.
$$\sum_{\g\in\G}\mu(\g) d(w,\g(w))<\infty\,.$$

The random walk $(Z_n)_n$ starting from the neutral element
$e$ associated with $\mu$ is defined by
the recursion relations:
$$Z_0=e\,;\, Z_{n+1}=Z_n\cdot X_{n+1}\, , $$
where $(X_n)$ is a sequence of independent and identically distributed
random variables of law $\mu$.
Thus, for each $n$, $Z_n$ is a random variable taking its values in $\G$.
We use the notation $Z_n(w)$ for the image of the base point $w\in X$ by $Z_n$.
The {\it rate of escape}, or {\it drift} of the random walk $Z_n(w)$ is the number $\ell$ defined as
$$\ell \fdefeq \lim_n \frac{d(w,Z_n(w))}n \,,$$
where the limit exists almost surely and in $L^1$ by the sub-additive ergodic Theorem
(J.\,Kingman \cite{kingman68}, Y.\,Derriennic \cite{derriennic75}).

If $\G$ is elementary, then its boundary is either empty or finite. In either case, there is no interest in looking
at properties at the boundary. We will assume from now on that $\G$ is non-elementary. In particular, $\G$ is non-amenable
so not only is the random walk always transient, $\ell$ is also positive (cf. \cite[\S\,7.3]{kai}).

There are different ways to prove that almost any trajectory of the random walk has a limit point
$Z_\infty(w)\in\partial X$. We recall below a theorem by V.\,Kaimanovich (cf. Theorem 7.3 in \cite{kai} and \S 7.4
therein) since it contains some information on the way
$(Z_n(w))$ actually tends to $Z_\infty(w)$ that will be used later.

\REFTHM{kai} {\em\bf (V.\,Kaimanovich)}.
Let $\Gamma$ be a non-elementary hyperbolic group and $(X,d)\in\DDD(\G)$, and
let us consider a symmetric
probability measure  $\mu$  with finite  first moment the support of which generates $\G$.
Then $(Z_n(w))$ almost surely converges to a point $Z_\infty(w)$ on the boundary.

For any $a\in\partial X$, we choose a quasigeodesic
$[w,a)$ from $w$ to $a$ in a measurable way.

For any $n$,
there is a measurable map $\pi_n$ from $\partial X$ to $X$ such that $\pi_n(a)\in [w,a)$,  and,
for almost
any trajectory of the random walk,
\begin{equation} \label{pin}
\lim_{n\to\infty}\frac{\vert Z_n(w)-\pi_n(Z_\infty(w))\vert }{n}=0\,.
 \end{equation}
\ENDTHM
The actual result was proved for geodesic metrics $d$.
Once proved in a locally finite Cayley graph, one may then use a quasi-isometry to
get the statement in this generality.

The estimate (\ref{pin}) will be improved in Corollary \ref{logcontrol}
under the condition that $d_G$ belongs to $\DDD(\G)$.

The  {\it harmonic measure} $\nu$ is then the law of $Z_\infty(w)$ i.e., it is
the probability measure on $\partial X$ such that
$\nu(A)$ is the probability that $Z_\infty(w)$ belongs to the set $A$.
More generally, we let $\nu_\g$ be the harmonic measure for the random walk started at
the point $\g(w)$, $\g\in\G$ i.e. the law of $\g ( Z_\infty(w))$. Comparing with the action of $\G$
on $\partial X$, we see that $\g^*\nu=\nu_{\g^{-1}}$.

\subsection{The Green metric} \label{greenmetric}

Let $\G$ be a countable group
and $\mu$ a symmetric law the support of which generates
$\G$.

For $x,y\in \G$, we define $F(x,y)$ as the probability
that a random walk starting from $x$ hits
$y$ in finite time i.e., the probability there is some $n$ such that $xZ_n=y$.
S.\,Blach\`ere and S.\,Brofferio \cite{bb} have defined the Green metric by
$$d_G(x,y) \fdefeq  -\log F(x,y)\,.$$

The Markov property implies that $F$ and the Green function $G$ satisfy 
$$ G(x,y)=F(x,y) G(y,y)\,.$$ 
Since $G(y,y)=G(e,e)$, we then get that 
$$F(x,y)=\frac{G(x,y)}{G(e,e)}$$ 
i.e. $F$ and $G$ only differ by a multiplicative contant and 
$$d_G(x,y)=\log G(e,e)-\log G(x,y)\,.$$

This function $d_G$ is known to be a left-invariant metric on $\G$
(see \cite{bb,bhm} for details).

We end this short introduction to the Green metric with the following folklore property.

\REFLEM{d_Gpropre} Let $\mu$ be a symmetric probability measure on $\G$ which defines
a transient random walk. Then $(\G,d_G)$ is a proper metric space i.e., balls of finite
radius are finite.\ENDLEM

\Proof It is enough to prove that $G(e,x)$ tends to $0$ as $x$ leaves any finite set.

Fix $n\ge 1$; by definition of convolution and by the Cauchy-Schwarz inequality,
$$
\mu^{2n}(x) = \sum_{y\in\G}\mu^{n}(y)\mu^{n}(y^{-1}x)\leq  \sqrt{\sum_{y\in\G}\mu^{n}(y)^2}\sqrt{\sum_{y\in\G}\mu^{n}(y^{-1}x)^2}
\,.$$

Since we are summing over the same set, it follows that
$$\sum_{y\in\G}\mu^{n}(y)^2=\sum_{y\in\G}\mu^{n}(y^{-1}x)^2$$
and the symmetry of $\mu$ implies that
$$
\sum_{y\in\G}\mu^{n}(y)^2 = \sum_{y\in\G}\mu^{n}(y)\mu^{n}(y^{-1}) = \mu^{2n}(e)\,.
$$

Therefore, $\mu^{2n}(x)\le \mu^{2n}(e)$.
Similarly,
$$
\mu^{2n+1}(x) = \sum_{y\in\G}\mu(y)\mu^{2n}(y^{-1}x)
\leq \sum_{y\in\G}\mu(y)\mu^{2n}(e)
\leq \mu^{2n}(e)\,.
$$

Since the walk is transient, it follows that $G(e,e)$ is finite, so, given $\ep>0$,
there is some $k\ge 1$ such that
$$\sum_{n\ge k} \mu^{2n}(e) \le \sum_{n\ge 2k} \mu^{n}(e)\le \ep\,.$$
On the other hand, since $\mu^n$ is a probability measure for all $n$, there is some finite
subset $K$ of $\G$ such that, for all $n\in\{0,\ldots,2k-1\}$, $\mu^n(K)\ge 1-\ep/(2k)$.
Therefore, if $x\not\in K$, then
$$
G(e,x) = \sum_{0\le n<2k}\mu^{n}(x)+  \sum_{n\ge 2k}\mu^{n}(x)
\le  \sum_{0\le n<2k}\mu^{n}(\G\setminus K)+  2\sum_{n\ge k }\mu^{2n}(e)
\le  \ep + 2\ep \, .
$$
The lemma follows.\endp

\subsection{The Martin boundary}\label{sec:martin}

Let $\G$ be a countable
group and $\mu$ be a symmetric probability measure
on $\G$. We assume that the support of $\mu$ generates $\G$
and that  the corresponding random walk is transient.

A non-negative function $h$ on $\G$ is {\it $\mu$-harmonic}
(harmonic for short) if, for all $x\in\G$,
$$h(x)=\sum_{y\in\G}h(y)\mu(x^{-1}y)\,.$$
A positive harmonic function $h$ is {\it minimal} if
any other positive harmonic function $v$ smaller
than $h$ is proportional to $h$.

The {\it Martin kernel} is defined  for all $(x,y)\in \Gamma \times \Gamma$ by
$$
K(x,y) \fdefeq \frac{G(x,y)}{G(e,y)} = \frac{F(x,y)}{F(e,y)} \, .
$$

We endow $\G$ with the discrete topology.
Let us briefly recall the construction of
the Martin boundary
$\partial_M \Gamma$:
let $\Psi:\Gamma\rightarrow C(\Gamma)$ be defined by
$y\longmapsto K_y=K(\cdot,y)$. Here $C(\Gamma)$ is the space of real-valued functions defined on $\Gamma$
endowed with the topology of pointwise convergence.
It turns out that $\Psi$ is injective and thus we may identify $\Gamma$ with its image.
The closure of $\Psi(\Gamma)$ is compact in $C(\Gamma)$ and, by definition,
$\partial_M \Gamma=\overline {\Psi(\Gamma) } \setminus \Psi(\Gamma)$ is the Martin boundary.
In the compact space $\Gamma \cup \partial_M \Gamma$, for any initial point $x$,
the random walk $Z_n(x)$ almost surely converges
to some random variable $Z_\infty(x)\in
\partial_M \Gamma$ (see for instance E.\,Dynkin \cite{dynkin69},
A.\,Ancona \cite{aa} or W.\,Woess \cite{wo}).

To every point $\xi\in\partial_M\G$ corresponds a positive
harmonic function $K_\xi$. Every minimal function arises in this way:
if $h$ is minimal, then there  are a constant $c>0$ and  $\xi\in\partial_M\G$ such
that $h=cK_\xi$. We denote by $\partial_m\G$ the
subset of $\partial_M\G$ consisting of (normalised) minimal positive harmonic functions.

Choquet's integral representation implies that, for any positive
harmonic function $h$, there is a unique probability measure
$\kappa^h$ on $\partial_m\G$ such that
$$h=\int K_{\xi}d\kappa^h(\xi)\,.$$

We will also use L.\,Na\"{\i}m's kernel $\Te$ on $\G\times\G$ defined
by $$\Te(x,y)\fdefeq\frac{G(x,y)}{G(e,x)G(e,y)}=\frac{K_y(x)}{G(e,x)}\,.$$
As the Martin kernel, Na\"{\i}m's kernel admits a continuous extension to $\G\times(\G\cup\partial_M\G)$.
In terms of the Green metric, one gets
\begin{equation}\label{lognaim}
\log\Te(x,y)= 2(x|y)^G_e -\log G(e,e)\,,
\end{equation}
where $(x|y)^G_e$ denotes the Gromov product with respect to the Green metric. 
See \cite{nai} for properties of this kernel.

We shall from now on assume that the Green metric $d_G$ is hyperbolic. Then
it has a visual boundary that we denote by $\partial_G \G$.
We may also compute the Busemann function in the metric $d_G$, say $\be_a^G$.
Sending $y$ to some point $a\in\partial_G \G$ in the equation
$ d_G(e,y)-d_G(x,y)=\log K(x,y) $, we get that
$\be_a^G(e,x)=\log K_a(x)$.

We now start
preparing the proof of \Th \ref{martin} in the next lemma and proposition.
We define an equivalence
relation $\sim_M$ on $\partial_M\G$: say that $\xi\sim_M\zeta$
if there exists a constant $C\ge 1$ such that
$$\frac{1}{C}\le\frac{K_\xi}{K_\zeta}\le C\,.$$
Given $\xi\in\partial_M\G$, we denote by $M(\xi)$ the class of $\xi$.

We first derive some properties of this equivalence relation:

\REFLEM{eqMprop}
\bit
\item[(i)] There exists a constant $E\ge 1$ such that for all sequences $(x_n)$ and $(y_n)$
in $\G$  converging to $\xi$ and $\zeta$ in $\partial_M\G$ respectively
 and such that $\Te(x_n,y_n)$ tends to infinity, then
$$\frac{1}{E}\le\frac{K_\xi}{K_\zeta}\le E\,;$$
in particular, $\xi\sim_M \zeta$.
\item [(ii)]For any $\xi\in\partial_M\G$, there is some $\zeta\in M(\xi)$
and a sequence $(y_n)$ in $\G$ which tends to some point $a\in\partial_G\G$ in the
sense of Gromov, to $\zeta\in\partial_M\G$ in the sense of Martin and such that $\Te(y_n,\xi)$
tends to infinity.
\item [(iii)]Let $\xi,\zeta\in\partial_M\G$. If $\zeta\notin M(\xi)$, then
there is a \nbhd $V(\zeta)$ of $\zeta$ in $\G$ and a constant $M$
such that $$K_\xi(x)\le M G(e,x)$$
for any $x\in V(\zeta)$.
\eit\ENDLEM

\Proof

\noindent (i) Fix $z\in\G$ and $n$ large enough so that $(x_n|y_n)^G_e\gg d_G(e,z)$;
we consider the approximate tree $T$ associated with $F=\{e,z,x_n,y_n\}$
and the $(1,C)$-quasi-isometry $\vp:(F,d_G)\to (T,d_T)$ (cf. \Th \ref{aptrees}).

On the tree $T$,we have
$$
|d_T(\vp(e),\vp(x_n))-d_T(\vp(z),\vp(x_n))|=|d_T(\vp(e),\vp(y_n))-d_T(\vp(z),\vp(y_n))|\,,
$$
so that $$ |(d_G(e,x_n)-d_G(z,x_n))-(d_G(e,y_n)-d_G(z,y_n))|\le 2C\,.$$
In terms of the Martin kernel,
$$|\log K_{x_n}(z)-\log K_{y_n}(z)|\le 2C\,.$$

Letting $n$ go to infinity yields the result.

\bigskip

\noindent (ii) Let $(y_n)$ be a sequence such that
$$\lim K_\xi(y_n)=\sup K_\xi\,.$$

Since $K_\xi$ is harmonic, the maximum principle implies that
$(y_n)$ leaves any compact set. But the walk is symmetric and transient so
Lemma \ref{d_Gpropre} implies that  $G(e,y_n)$ tends to $0$.

Furthermore, for $n$ large enough, $K_\xi(y_n)\ge K_\xi(e)=1$,
so that
$$\Te(y_n,\xi)\ge \frac{1}{G(e,y_n)}\to\infty\,.$$

Let $(x_n)$ be a sequence in $\G$ which tends to $\xi$.
For any $n$, there is some $m$ such that
$$|K_\xi(y_n)-K_{x_m}(y_n)|\le G(e,y_n)\,.$$
It follows that
$$\Te(y_n,x_m)\ge \Te(y_n,\xi) - \frac{|K_\xi(y_n)-K_{x_m}(y_n)|}{G(e,y_n)}\ge \Te(y_n,\xi) -1\,.$$

Therefore, applying part (i) of the lemma,
we see that  any limit point of $(y_n)$  in $\partial_M\G$ belongs to $M(\xi)$.

Moreover, for any such limit point $\zeta\in\partial_M\G$, we get that
$$\Te(y_n,\zeta)\ge \frac{1}{E}\Te(y_n,\xi)\,.$$
Applying the same argument as above,
we see that, for any $M >0$, there is some $n$ and $m_n$
such that, if $m\ge m_n$ then
$$\Te(y_n,y_m)\ge M -1\,.$$
From (\ref{lognaim}) we conclude, using a diagonal procedure, that there exist a subsequence $(n_k)$
such that $(y_{n_k})$ tends to infinity in the Gromov topology.

\bigskip

\noindent (iii) Since $\zeta\notin M(\xi)$, there is a \nbhd $V(\zeta)$
and a constant $M$ such that $\Te(x,\xi)\le M$ for all $x\in V(\zeta)$.
Otherwise, we would find $y_n\to\zeta$ with $\Te(y_n,\xi)$ going
to infinity, and the argument above would imply $\zeta\in M(\xi)$.
Therefore, $$K_\xi(x)\le M G(e,x)\,.$$\endp

\REFPROP{minimal} Every Martin point is minimal.
\ENDPROP

\Proof We observe that if $K_\xi$ is minimal, then
$M(\xi)=\{\xi\}$. Indeed, if $\zeta\in M(\xi)$, then
$$K_\xi\ge K_\xi -\frac{1}{C} K_\zeta\ge 0$$
for some constant $C\ge 1$. The minimality of $K_\xi$
implies that  $K_\xi$ and $K_\zeta$ are proportional and, since their value at $e$ is $1$,
$K_\xi=K_\zeta$ i.e., $\xi=\zeta$.

Let $\xi\in\partial_M\G$. There is a unique probability measure
$\kappa^\xi$ on $\partial_m\G$ such that
$$K_\xi=\int K_\zeta d\kappa^\xi(\zeta)\,.$$
By Fatou-Doob-Na\"{\i}m Theorem, for $\kappa^\xi$-almost every
$\zeta$, the ratio $G(e,x)/K_\xi(x)$ tends to $0$ as $x$ tends
to $\zeta$ in the fine topology \cite[Thm. II.1.8]{aa}. From Lemma \ref{eqMprop} (iii),
it follows that $\kappa^\xi$ is supported by $M(\xi)$. In particular,
$M(\xi)$ contains a minimal point.\endp

{\noindent\sc Proof of \Th \ref{martin}.} Since every Martin point is minimal,
Lemma \ref{eqMprop}, (ii), implies that for every $\xi\in\partial_M\G$,
there is some sequence $(x_n)$ in $\G$ which tends to $\xi$ in the Martin
topology and to some point $a$ in the hyperbolic boundary as well.

Let us prove that the point $a$ does not depend on the sequence.
If $(y_n)$ is another sequence tending to $\xi$, then
 $$\limsup_{n,m\to\infty}\Te(x_n,y_m)=\infty$$
because $\Te(\xi,x_n)$ tends to infinity. Therefore, there is
a subsequence of $(y_n)$ which tends to $a$ in the Gromov topology.
Since we have only one accumulation point, it follows that
$a$ is well-defined. This defines a map $\phi:\partial_M\G\to\partial_G\G$.

Now, if $(x_n)$ tends to $a$ in the Gromov topology,
then it has only one accumulation point in the Martin boundary as well by Lemma \ref{eqMprop}, (i). 
So the map $\phi$ is injective.
The surjectivity follows from the compactness of $\partial_M\G$.

To  conclude the proof, it is enough to prove the continuity of $\phi$
since $\partial_M\G$ is compact. Let $M>0$ and $\xi\in\partial_M\G$
be given.
We consider a sequence $(x_n)$ which tends to $\xi$
as in Lemma \ref{eqMprop}.
Let $C$ be the constant given
by \Th \ref{aptrees} for $4$ points. We pick $n$ large enough so
that $(x_n|\phi(\xi))^G_e\ge M+2C+\log 2$. Let $$A=\min\{K_\xi(x), x\in B_G(e, d_G(x_n,e))\}.$$

Let $\zeta\in\partial_M\G$ such that $|K_\xi-K_\zeta|\le (A/2)$ on
$B_G(e, d_G(x_n,e))$. It follows that $$1/2\le \frac{K_\zeta}{K_\xi}\le 3/2\,.$$

Approximating $\{e,x_n,\phi(\xi),\phi(\zeta)\}$ by a tree,
we conclude that $(\phi(\xi)|\phi(\zeta))^G_e\ge M$, proving the continuity
of $\phi$.\endp

\subsection{Hyperbolicity of the Green metric}

We start with a characterisation of the hyperbolicity of the Green metric
in the quasi-isometry class of the group.

\REFPROP{d_GdansD(G)}  Let $\G$ be a non-elementary hyperbolic group and $\mu$ a symmetric probability
measure with Green function $G$. We fix a finite generating set
$S$ and consider the associated word metric $d_w$.
The Green metric $d_G$ is quasi-isometric to $d_w$ and hyperbolic if and only if the following
two 
conditions are satisfied.
\bit
\item[(ED)] There are positive constants $C_1$ and  $c_1$ such that, for all $\g\in\G$,
$$G(e,\g)\le C_1 e^{-c_1 d_w(e,\g)}$$
\item[(QR)] For any $r\ge 0$, there exists a positive constant $C(r)$ such that
$$G(e,\g)\le C(r) G(e,\g')G(\g',\g)$$
whenever $\g,\g'\in\G$ and $\g'$ is at distance at most $r$ from a  $d_w$-geodesic segment between $e$ and $\g$.
\eit
\ENDPROP

{\noindent\bf Remark.}
Even though hyperbolicity is an invariant property under quasi-isometries between geodesic metric spaces,
this is not
the case when we do not assume the spaces to be geodesic (see the appendix).

\bigskip

\Proof We first assume that $d_G\in\DDD(\G)$. The quasi-isometry property implies that condition (ED) holds.
The second condition (QR) follows from \Th \ref{main}.

Indeed, since $d_G$ is hyperbolic and quasi-isometric to a word distance, then
$(\G,d_G)$ is quasiruled. This is sufficient to ensure that condition (QR) holds for
$r=0$.   The general case $r\geq 0$ follows: let $y$ be the closest point to $\g'$ on a geodesic
between $e$ and $\g$ and note that
$$ d_G(e,\g')+d_G(\g',\g)
\le d_G(e,y)+d_G(y,\g)+2d_G(y,\g')
\le
\log C(0)+d_G(e,\g)+2d_G(y,\g')\,.$$
Thus one may choose
$C(r)=C(0) \exp(2c)$ where $c=\sup d_G(y,\g')$ for all pair $y,\g'$ at distance less than $r$.
This last $\sup$ is finite because $d_G$ is quasi-isometric to a word metric.

\bigskip

For  the converse, we assume that
both conditions (ED) and (QR) hold and
let \\ $C=\max\{d_G(e,s),\ s\in S\}$. For any $\g\in\G$, we consider a
$d_w$-geodesic $\{\g_j\}$ joining $e$ to $\g$.
It follows that $$d_G(e,\g)\le \sum_j d_G(\g_j,\g_{j+1})\le Cd_w(e,\g)\,.$$
From (ED), we obtain $$d_G(e,\g)\ge  c_1 d_w(e,\g)-\log C_1\,.$$
Since both metrics are left-invariant, it follows that $d_w$ and $d_G$ are quasi-isometric.

\bigskip

Condition (QR) implies
that $d_w$-geodesics are not only quasigeodesics for $d_G$, but also
quasirulers,  cf. Appendix \ref{app}. Indeed, since the two functions $F$ and $G$ only
differ by a multiplicative factor, condition (QR) implies
that there is a constant $\tau$ such that,
for any $d_w$-geodesic segment $[\g_1,\g_2]$
and any $\g\in [\g_1,\g_2]$,
 we have
$$ d_G(\g_1,\g)+d_G(\g,\g_2)\le 2\tau + d_G(\g_1,\g_2)\,.$$
\Th \ref{main}, (iii) implies (i), implies that $(\G,d_G)$ is a hyperbolic space.\endp

To prove the first statement of \Th \ref{main1}, it is now enough to establish the following lemma.

\REFLEM{d_Gqisom} Let $\G$ be a non-elementary hyperbolic group, and $\mu$ a symmetric probability
measure with finite exponential moment. Then condition (ED) holds.
\ENDLEM

When $\mu$ is finitely supported, the lemma was proved by S.\,Blach\`ere and S.\,Brofferio
using the Carne-Varopoulos estimate \cite{bb}.

\bigskip

\Proof Let us fix a word metric $d_w$ induced by a finite generating set $S$, so that $d_w\in\DDD(\G)$.

Since $\G$ is non-amenable,
Kesten's criterion implies that there are positive constants $C$ and $a$ such that
\begin{equation} \label{kesten}
\forall \gamma \in \G , \;\;\; \mu^n(\gamma) \leq \mu^n (e) \leq C e^{-an} \, .
\end{equation}
For a proof, see \cite[Cor. 12.5]{wo}.

\bigskip

We assume that $\EE [\exp \la d_w(e,Z_1)]=E<\infty$ for a given $\la>0$.
For any $b>0$, it follows from the exponential Tchebychev inequality
that $$\PP\left[\sup_{1\le k\le n} d_w(e,Z_k)\ge nb\right]\le e^{-\la b n}\EE\left[ \exp\left(\la \sup_{1\le k\le n} d_w(e,Z_k)\right)\right]\,.$$
But then, for $k\le n$,  $$ d_w(e,Z_k) \le \sum_{1\le j\le n-1} d_w(Z_j,Z_{j+1})
= \sum_{1\le j\le n-1} d_w(e,Z_j^{-1}Z_{j+1})\,.
$$

The increments $(Z_j^{-1}Z_{j+1})$ are independent random variables
and all follow the same law as $Z_1$. Therefore
\begin{equation} \label{increment}
\PP\left[\sup_{1\le k\le n} d_w(e,Z_k)\ge nb\right] \le e^{-\la b n} E^n= e^{(-\la b+\log E)n}\,.
\end{equation}
We choose $b$ large enough so that $c\fdefeq -\la b+\log E <0$.

 We have
$$G(e,\g)=\sum_n\mu^n(\g)=\sum_{1\le k\le |\g|/b}\mu^k(\g) +\sum_{ k> |\g|/b}\mu^k(\g)\,,$$
where we have set $|\g|=d_w(e,\g).$
The estimates  (\ref{increment}) and (\ref{kesten}) respectively imply that
\begin{eqnarray*}
\sum_{1\le k\le |\g|/b}\mu^k(\g) & \le & \frac{|\g|}{b} \sup_{1\le k\le |\g|/b} \mu^k(\g) \le  \frac{|\g|}{b}
\PP [\exists k\le |\g|/b \mbox{ s.t. } Z_k=\g]\\
 & \le & \frac{|\g|}{b} \PP \left[ \sup_{1\le k\le |\g|/b} d_w(e,Z_k) \ge |\g|\right] \lesssim |\g|e^{-c|\g|}
\end{eqnarray*}
and
$$\sum_{ k> |\g|/b}\mu^k(\g)\lesssim e^{-(a/b)|\g|}\,.$$
Therefore, (ED) holds.\endp

When $\G$ is hyperbolic and $\mu$ has finite support, A.\,Ancona \cite{aa} proved that
the Martin boundary is homeomorphic to the visual boundary $\partial X$. The key point
in his proof is the following estimate (see \cite[Thm. 27.12]{wo} and \Th \ref{martin}).

\REFTHM{anc} {\em\bf (A.\,Ancona)}
Let $\Gamma$ be a non-elementary hyperbolic group,
$X$ a locally finite Cayley graph endowed with a geodesic metric $d$ so that
$\G$ acts canonically by isometries, and let $\mu$ be a finitely supported symmetric probability measure
the support of which generates $\Gamma$.
For any $r\ge 0$, there is a constant $C(r)\ge 1$ such that $$F(x,v)F(v,y)\le F(x,y)\le C(r)F(x,v)F(v,y)$$
whenever $x,y\in X$ and $v$ is at distance at most $r$ from a geodesic segment between $x$ and $y$.\ENDTHM

This implies together with Lemma \ref{d_Gqisom} that when $\mu$ is finitely supported, both conditions
(ED) and (QR) hold. Therefore, Proposition \ref{d_GdansD(G)} implies that $d_G\in\DDD(\G)$.
We have just established the first statement of Corollary \ref{cmain1}.

\subsection{Martin kernel vs Busemann function: end of the proof of Theorem \ref{main1}}
We assume that $X=\G$ equipped with the Green metric $d_G$ belongs to $\DDD(\G)$ throughout this paragraph.

{\noindent\bf Notation.} When we consider notions with respect to $d_G$,
we will add the exponent $G$ to distinguish them from the same
notions in the initial metric $d$.  Thus Busemann functions for $d_G$ will be written $\be_a^G$.
The visual metric on $\partial X$ seen from $w$ for the original metric will be denote by  $d_\ep$, and by $d_\ep^G$ for the
one coming from $d_G$. Balls at infinity will  be denoted by $B_\ep$ and $B_\ep^G$.

\bigskip

Let us recall that the Martin kernel is defined by $$K(x,y)=\frac{F(x,y)}{F(w,y)}
 =\exp\left\{ d_G(w,y)- d_G(x,y)\right\}. $$

By definition of the Martin boundary $\partial_M X$, the kernel $K(x,y)$ continuously extends
to a $\mu$-harmonic positive function $K_a(\cdot)$
when $y$ tends to a point $a\in \partial_M X$.
We recall that,  by  \Th \ref{martin}, we may - and will -
identify  $\partial_M X$ with the visual boundary $\partial X$.

As we already mentioned $\Gamma$ acts on $\partial_M X$,
so on its harmonic measure and we have
$\g^*\nu=\nu_{\g^{-1}}$.
Besides,
see e.g. G.\,Hunt \cite{hunt60} or W.\,Woess \cite[Th. 24.10]{wo} for what follows,
$\nu$ and $\nu_\g$ are absolutely continuous and their
Radon-Nikodym derivatives satisfy
$$\frac{d\nu_\g}{d\nu}(a) =  K_a(\g(w))\,.$$

We already computed the Busemann function in the metric $d_G$
 in part \ref{sec:martin}:
$\be_a^G(w,x)=\log K_a(x)$. Thus we have proved that
$$ \frac{d\g^*\nu}{d\nu}(a) = \exp \be_a^G(w,\g^{-1}w)\, . $$

It follows at once that $\nu$ is a \qc measure on $(\partial X, d_{\ep}^G)$ of dimension $1/\ep$.
Actually, $\nu$ is even a {\it conformal measure} since we have a genuine equality above.
Therefore $\nu$ belongs to the Patterson-Sullivan class associated with the metric $d_G$.
According to \Th \ref{reg}, it is in particular comparable to the Hausdorff measure for the corresponding visual metric.
This ends both the proofs of \Th \ref{main1} and of Corollary \ref{cmain1}.

We note that, comparing the statements in  \Th \ref{main1} (ii) and \Th \ref{reg}, we
recover the equality  $v_G=1$ already noticed in \cite{bb} for random walks on non-amenable groups.
See also \cite{bhm}.

\subsection{Consequences} We now draw consequences  of the hyperbolicity of the Green metric.

We refer to the appendices for properties of quasiruled spaces.

\subsubsection{Deviation inequalities}

We study the lateral deviation of the position
of the random walk with respect to the quasiruler $[w,Z_{\infty}(w))$ where,
for any $x \in X$ and $a \in \partial X$,  we chose an arbitrary quasiruler $[x,a)$ from $x$ to $a$
in a measurable way.

\REFPROP{control} Assume that $\G$ is a non-elementary hyperbolic group, $(X,d)\in\DDD(\G)$, and  $\mu$ is a symmetric law so that the associated
Green metric belongs to $\DDD(\G)$.
The following holds
\bit
\item[(i)] There is a positive constant $b$ so that,
for any $D\ge 0$ and $n\ge 0$,
$$ \PP[d(Z_n(w), [w,Z_\infty(w)))\ge D] \lesssim e^{-bD}\,.$$
\item[(ii)] There is a constant $\tau_0$ such that
for any positive integers $m,n,k$,
$$\EE[(Z_m(w)|Z_{m+n+k}(w))_{Z_{m+n}(w)}]\le \tau_0\,.$$\eit\ENDPROP

\Proof

\noindent Proof of (i).
Observe that 
\begin{eqnarray*}
\PP[d(Z_n(w), [w,Z_\infty(w)))\ge D] & = & \sum_{z\in X} \PP[d(Z_n(w), [w,Z_\infty(w)))\ge D \, ,\, Z_n(w)=z]\\
& = & \sum_{z\in X} \PP[d(z,[w,Z_n^{-1}Z_\infty(z)))\ge D \, ,\, Z_n(w)=z]\\
& = & \sum_{z\in X} \PP[d(z,[w,Z_n^{-1}Z_\infty(z)))\ge D ] \PP[ Z_n(w)=z]\\
& = & \sum_{z\in X} \PP[d(z,[w,Z_\infty(z)))\ge D ] \PP[ Z_n(w)=z]
\end{eqnarray*}

The second equality holds because $\gamma w=z$ implies that $\gamma^{-1} Z_\infty(z)=Z_\infty(w)$. 
The third equality comes
from the independence of $Z_n=X_1X_2\cdots X_n$ and $Z_n^{-1}Z_\infty =X_{n+1}X_{n+2}\cdots$.
The last equality uses the fact that $Z_n^{-1}Z_\infty$ and $Z_\infty$ have the same law.

On the event $\{d(z,[w,Z_\infty(z)))\geq D\}$,  we have in particular $d(w,z)\geq D$ and
we can pick $x \in [w,z)$ such that $d(z,x)=D +O(1)$.
Then, because the triangle $(w,z,Z_{\infty}(z))$ is thin and since $d(z,[w,Z_\infty(z)))\geq D$,
we must have $Z_\infty(z) \in \mho_z(x,R)$. As usual $R$ is a constant that does not depend
on $z$, $D$ or $Z_{\infty}(z)$.
We now apply the lemma of the shadow Lemma \ref{ashadow} to the Green metric to deduce that

$$
\PP[d(z,[w,Z_\infty(z)))\ge D ] \leq \PP^z[Z_\infty(z)\in\mho_z(x,R)]=
\nu_z(\mho_z(x,R)) \lesssim e^{-d_G(z,x)}\,.$$

Finally, using the quasi-isometry between $d$ and $d_G$,
it follows that
$$ \PP[d(Z_n(w), [w,Z_\infty(w)))\ge D] \lesssim e^{-bD}\,.$$

\noindent Proof of (ii).
Using the independence of the increments of the walk, one may
first assume that $m=0$.

Let us  choose $Y_n(w)\in [w,Z_{\infty}(w))$ such that
$d(w,Y_n(w))$ is as close from
$(Z_n(w)|Z_{\infty}(w))$ as possible. Since the space $(X,d)$ is quasiruled, it follows that
$d(w,Y_n(w))=(Z_n(w)|Z_{\infty}(w)) +O(1)$.

(We only use Landau's notation $O(1)$ for estimates that are uniform
with respect to the trajectory of $(Z_n)$. Thus the line just above means that there exists a
\underline{deterministic} constant $C$ such that
$$|d(w,Y_n(w)) - (Z_n(w)|Z_{\infty}(w))|\le C\,.$$ The same convention applies to the rest of the proof.)

Let us define
$$A_0=\{ d(w,Y_n(w))\le d(w,Y_{n+k}(w))\}$$
and, for $j\ge 1$,
$$A_j= \{j-1 <  d(w,Y_n(w))- d(w,Y_{n+k}(w))\le j\}\,.$$
Approximating $\{w,Z_n(w),Z_{n+k}(w),Z_\infty(w)\}$ by a tree, it follows
that, on the event $A_0$,
$$(w|Z_{n+k}(w))_{Z_{n}(w)}\le
d(Z_n(w),[w,Z_\infty(w)))+O(1)$$
and that, on the event $A_j$,
$$(w|Z_{n+k}(w))_{Z_{n}(w)}\le
d(Z_n(w),[w,Z_\infty(w)))+j+O(1)\, .$$

Therefore
$$\EE[(w|Z_{n+k}(w))_{Z_{n}(w)}]\leq  \EE [d(Z_n(w),[w,Z_\infty(w)))] +\sum_{j\ge 1}j\PP(A_j)+O(1)\,.$$

If $d(w,Y_n(w))- d(w,Y_{n+k}(w))\ge j$ then $d(Z_{n+k}(w),[Z_n(w),Z_\infty(w)))\ge j$
so that $$\PP(A_{j+1})\le \PP[d(Z_{n+k}(w),[Z_n(w),Z_\infty(w)))\ge j]\,.$$

Using (i) for the random walk starting at $Z_n(w)$, we get

$$\sum_{j\ge 1}j\PP(A_j)\lesssim 1\,.$$

On the other hand,
$$\EE [d(Z_n(w),[w,Z_\infty(w)))]
=\int_0^{\infty} \PP[d(Z_n(w), [w,Z_\infty(w)))\ge D] \, dD \lesssim \int_0^{\infty} e^{-bD} \, dD =1/b \, .
$$
The proposition follows.
\endp

We now improve
the estimate (\ref{pin}) in \Th \ref{kai} when $d_G \in \DDD(\G)$.

\REFCOR{logcontrol}  Let $\G$ be a non-elementary hyperbolic group, $(X,d)\in\DDD(\G)$ and $\mu$ a symmetric law such that $d_G\in\DDD(\G)$, then we have
\begin{equation} \label{latdev}
\limsup \frac{d(Z_n(w),[w,Z_\infty(w)))}{\log n}<\infty \;\;\; \PP \mbox{ a.s.}
\end{equation}
\ENDCOR

\Proof It follows from Proposition \ref{control} that we may find a constant $\kappa >0$ so that
$$\PP[d(Z_n(w),[w,Z_\infty(w)))\ge \kappa \log n]\le \frac{1}{n^2}\,.$$
Therefore, the Borel-Cantelli lemma implies that
$$\limsup \frac{d(Z_n(w),[w,Z_\infty(w)))}{\log n}<\infty \;\;\; \PP \mbox{ a.s.}$$
and the corollary follows.\endp

{\noindent\bf Remark.} After a first version of this paper was publicised, M.\,Bjorklund also used 
the hyperbolicity of the Green metric to 
prove a Central Limit Theorem for $d_G(w,Z_n(w))$, see \cite{bjork}.

\subsubsection{Escape of the random walk from balls}
We assume here that $\mu$ is a symmetric and finitely supported probability measure on
a non-elementary hyperbolic group $\G$ and that the support of $\mu $ generates $\G$.
We want to compare the harmonic measure with the uniform measure
on the spheres for the Green metric.
We define the (exterior) sphere of the ball $B_G(w,R)$ by

$$
\partial B_G(w,R) \fdefeq \{x\in X \; :\; x\not\in B_G(w,R) \mbox{ and }
\exists \gamma \in Supp(\mu) \mbox{ s.t. } \gamma^{-1} (x)\in B_G(w,R)\} \, .
$$
The harmonic measure $\nu_R$ on $\partial B_G(w,R)$ is the law of the first point visited outside $B_G(w,R)$.

As the volume of the sphere $\partial B_G(w,R)$  equals $e^R$ up to a multiplicative constant (see \cite{bb}),
we need to compare $\nu_R (\cdot)$ with $e^{-R}$.
In other words, 
we have to bound the ratio between
the measure $\nu_R(\cdot)$ and the hitting probability $F(w,\cdot)$.
 Observe that, in principle, there could be points on the sphere that are visited
by the walk a long time after it  left the ball.
We shall see that this scenario can only take place on a finite scale.

In the following  we only consider quasigeodesics for $(X,d)$ and $(X,d_G)$
that are geodesics for a given word metric $d_w\in\DDD(\G)$.

\REFPROP{compmesharm} There exist positive constants $C_1<1$ and $C_2$ such that for any positive real $R$, the harmonic measure $\nu_R$ on the sphere
$\partial B_G (w,R)$ satisfies
$$
\forall x \in \partial B_G (w,R), \; \exists y \in B_G(x,C_1) \cap \partial B_G (w,R) \;\mbox{ s.t.} \; C_2 e^{-R} \leq \nu_R(y) \leq
e^{-R} \, .
$$
\ENDPROP

\Proof The upper bound (valid for any $x \in \partial B_G (w,R)$)
obviously follows from the definition of the Green metric:
if $y \not\in B_G(w,R)$,
then
$$
\nu_R (y)\leq F(w,y) = \exp (-d_G(w,y)) \leq e^{-R}\, .
$$

For the lower bound, we consider a quasigeodesic from $w$ to $x$ and denote by $y$
the first point of $\partial B_G (w,R)$ along that path.
Since  $\mu$ has finite support, $d_G(w,x)$ and $d_G(w,y)$ {only differ by} an additive constant.
The quasiruler property then implies that $y$ is at a bounded distance from $x$.

Let  $\EEE=\EEE(R)$ denote the set of points $z\in \partial B_G (w,R)$ such that
there is a quasigeodesic reaching $z$ from $w$ entirely
contained in $B_G(w,R)$ (except for the last step toward $z$).

Let $z\in \EEE$.
Since $y$ and $z$ belong to $\partial B_G (w,R)$, 
then $d_G(w,z)$ and $d_G(w,y)$ only differ by an additive constant and 
we have
\begin{equation} \label{2k}
d_G(y,z) \geq   d_G(y,z)+ (d_G(w,y) - d_G (w,z) -C)= 2(w|z)_y-C
\end{equation}

Let $k_0$ be an integer and define
$$
\EEE_0\fdefeq \{z\in \EEE \; : \; (w|z)_y \leq k_0\}
$$
and for all integer $k \geq k_0$,
$$
\EEE_k\fdefeq \{z\in \EEE \; : \; k<(w|z)_y \leq k+1\} \, .
$$
We denote by $\tau_R$ the first hitting time of $\partial B_G(w,R)$ by the random walk and by $\tau_y$
the first hitting time of $y$.
Then
$$
F(w,y) = \PP [ \tau_y <\infty, \, Z_{\tau_R}(w) \in \EEE_0] + \sum_{k=k_0}^{\infty} \sum_{z\in \EEE_k}
\PP [ \tau_y <\infty ,\, Z_{\tau_R}(w) =z]
$$
At this point, we need to use the Strong Markov property to say that 
once we know that $Z_{\tau_R}(w) =z$ and $z \neq y$,
the hitting time of $y$ must occur after $\tau_R$. Then,
the finiteness of $\tau_y$ depends only on the position $z$ disregarding the behavior of
the random walk up to time $\tau_R$.
Namely,
$$
\PP [ \tau_y <\infty ,\, Z_{\tau_R}(w) =z]=\PP^z [ \tau_y <\infty ]\PP [  Z_{\tau_R}(w) =z] \, .
$$
Using (\ref{2k}), the definition of $(\EEE_k)$ and  the inequality
$\PP [  Z_{\tau_R}(w) =z]\le \PP [\tau_z <\infty] \leq e^{-R}$, we get that
\begin{equation} \label{sum2k}
F(w,y) \leq \PP [ Z_{\tau_R}(w) \in \EEE_0] + C \sum_{k=k_0}^{\infty} e^{-2k} e^{-R} \# \EEE_k \, .
\end{equation}
We need an upper bound on $\# \EEE_k$.
Take $z \in \EEE_k$, and
let $y_{R-k}$ be the point at distance $R-k$ from $w$ along the quasigeodesic $[w,y]$.

As the triangle $(w,z,y)$ is thin, the center of the associated approximate tree
is at a bounded distance from the point $y_{R-k}$.
Then, since for any $z$ in $\EEE_k$, $(w|y)_z-k$ is bounded by a constant,
the set $\EEE_k$ is  therefore included in the ball $B_G(y_{R-k},k+C)$ for some constant $C$.
Thus $\# \EEE_k \lesssim e^k$ and
\begin{equation} \label{sum2k0}
C \sum_{k=k_0}^{\infty} e^{-2k} e^{-R} \# \EEE_k \leq C(k_0) e^{-R}
\end{equation}
with $C(k_0)$ tending to $0$ when $k_0$ tends to infinity.

As $\mu$ is finitely supported, $\partial B_G(w,R)$ is at a bounded distance from $B_G(w,R)$.
So $y \in B_G(w,R+C(\mu))$ and $F(w,y) \geq e^{-C(\mu)} e^{-R}$.
Now choose $k_0$ so that $C(k_0) < (1/2) e^{-C(\mu)}$ and take $R>k_0$.
Then (\ref{sum2k}) and (\ref{sum2k0}) give
\begin{equation} \label{harmeps0}
\PP [ Z_{\tau_R}(w) \in \EEE_0] \geq \frac12 e^{-C(\mu)} e^{-R}
\, .
\end{equation}
We conclude that $\nu_R(\EEE_0)\gtrsim e^{-R}$.
Take $y'\in \EEE_0$ so that $(w|y')_y\leq k_0$.
By the definition of the set $\EEE$ and by the thinness of the triangle $(w,y,y')$,
there exists a path joining $y$ and $y'$ within $B_G(w,R)$ of length at most $c(k_0)$,
 a constant depending only on $k_0$ and $\delta$.
Therefore, there exists a constant $c'(k_0,\mu)$ such that
$$
\nu_R(y)\geq \nu_R(y')c'(k_0,\mu) \, .
$$
Finally, as $\#\EEE_0$ is bounded above by a constant, (\ref{harmeps0}) gives
$$
\nu_R(y) \gtrsim \sum_{y'\in \EEE_0} \nu_R(y')=\nu_R(\EEE_0) \gtrsim \ep^{-R} \, .
$$
\endp

\noindent {\bf Remark.}
Proposition \ref{compmesharm} says that the harmonic measure on spheres
is well spread out and that the harmonic measure of a bounded domain
of the sphere of radius $R$ if $e^{-R}$ up to a multiplicative constant.
Approximating the balls of $\partial X$ by shadows, we get that $\nu$ is
Alfhors-regular of dimension $1/\ep$,  hence quasiconformal.
Therefore, we get an alternative proof of the second statement of Theorem \ref{main1}
when $\mu$ has finite support.

\subsubsection{The doubling condition for the harmonic measure}

Let us recall that a measure $m$ is said to be {\it doubling} if there exists a constant $C>0$ such that,
for any ball $B$
of radius at most the diameter of the space 
then $m(2B)\le Cm(B)$.

\REFPROP{cdoub} Let $\G$ be a non-elementary hyperbolic group, $(X,d)\in\DDD(\G)$
and let $\mu$ be a symmetric law such that
$d_G\in\DDD(\G)$.
The harmonic measure is doubling with respect to the visual measure
$d_{\ep}$ on $\partial X$.\ENDPROP

\Proof The modern formulation of Efremovich and Tichonirova's \th (cf. \Th 6.5 in \cite{bsc}
and references therein) states that quasi-isometries between hyperbolic proper geodesic
spaces $\Phi:X\to Y$ extend as quasisymmetric maps $\phi:\partial X\to \partial Y$
between their visual boundaries i.e., there is an increasing \homeo $\eta:\R_+\to\R_+$ such that
$$|\phi(a)-\phi(b)|\le \eta(t) |\phi(a)-\phi(c)|$$ whenever $|a-b|\le t|a-c|$.

Since $d_G\in\DDD(\G)$, the spaces involved are visual. Thus, the statement remains true since
we may still approximate properly the space by trees, cf. Appendix \ref{appsha}.

Since $(X,d)$ and $(X,d_G)$ are quasi-isometric, the boundaries  are thus quasisymmetric with respect
to $d_\ep$ and $d_\ep^G$. Furthermore, $\nu$ is doubling with respect to $d_\ep^G$ since it is Ahlfors-regular,
and this property is preserved under quasisymmetry.\endp

Basic properties on quasisymmetric maps include \cite{hei}. More information on boundaries of hyperbolic
groups, and the relationships between hyperbolic geometry and conformal geometry can be found in
\cite{bp,bk}.


\section{Dimension of the harmonic measure on the boundary of a hyperbolic metric space}

\Th \ref{main2} will follow from Proposition \ref{dimp} and Proposition \ref{dimH}.

We recall the definition of the rates of escape $\ell$ and $\ell_G$ of the random walk with respect to $d$ or $d_G$ respectively.
$$\ell \fdefeq \lim_n \frac{d(w,Z_n(w))}n \mbox{ and } \ell_G \fdefeq \lim_n \frac{d_G(w,Z_n(w))}n \, .$$

We will first prove

\REFPROP{dimp}
Let $\G$ be a non-elementary hyperbolic group and let $(X,d)\in\DDD(\G)$.
Let $\mu$ be a  symmetric probability measure on $\Gamma$ the support of which generates $\Gamma$ such that
$d_G\in\DDD(\G)$ and with finite first moment
$$\sum_{\g\in\G} d_G(w,\g(w))\mu(\g)<\infty\,.$$
Let $\nu$ be the harmonic measure seen from $w$ on $\partial X$.

For $\nu$-a.e. $a\in\partial X$,
$$\lim_{r\to 0} \frac{ \log \nu (B_{\ep}(a,r))}{\log r} = \frac{\ell_G}{\ep\ell}\,,$$
where $B_{\ep}$ denotes the ball on $\partial X$ for the visual metric $d_{\ep}$.
\ENDPROP

\noindent {\bf Remark.}
Recall from \cite{bhm} that $\mu$ having finite first moment with respect to the Green metric
is a consequence 
of $\mu$ having finite entropy.

\Proof
It is convenient to introduce an auxiliary word metric $d_w$ which is of course geodesic.
We may then consider the visual quasiruling structure $\GGG$ induced by the $d_w$-geodesics
for both metrics $d$ and $d_G$ via the identity map, cf. the appendix.

We combine Propositions \ref{shadowVSball} and \ref{shadowVSballbis} 
to get that, for  a fixed but large enough $R$, for any $a \in \partial X$ and $x\in [w,a)\subset\GGG$
$$
B_{\ep}(a,(1/C) e^{-\ep d(w,x)}) \subset \mho_{\GGG} (x,R) \subset B_{\ep}(a,C e^{-\ep d(w,x)})
$$
and
$$
B_{\ep}^G(a,(1/C) e^{-\ep d_G(w,x)}) \subset\mho_{\GGG} (x,R) \subset B_{\ep}^G(a,C e^{-\ep d_G(w,x)})
$$
for some positive constant $C$.
We recall that the shadows $\mho_{\GGG} (x,R)$ are defined using geodesics for the word metric $d_w$.

The doubling property of $\nu$ with respect to the visual metric $d_\ep$  implies that
\begin{equation} \label{compballs}
\nu (B_{\ep}(a,C e^{-\ep d(w,x)})) \asymp \nu (\mho_{\GGG}(x,R))
\end{equation}
for any $x \in [w,a)$.

Let $\eta>0$; by definition of the drift, there is a set of full measure with respect to
the law of the trajectories of the random walk, in which for any sequence $(Z_n(w))$ and for $n$
large enough, we have $|d(w,Z_n(w))-\ell n|\le \eta n$ and
$|d_G(w,Z_n(w))-\ell_G n|\le \eta n$.

From Theorem \ref{kai} applied to the metrics $d$  and $d_G$,  we get that, for $n$ large enough,
$d(Z_n(w),\pi_n(Z_\infty(w)))\le \eta n$ and $d_G(Z_n(w),\pi_n(Z_\infty(w))) \le \eta n$.

We conclude that

\begin{equation}\label{disdri}
\left\{\begin{array}{l} |d(w,\pi_n(Z_\infty(w)))-\ell n|\le 2 \eta n\\
|d_G(w,\pi_n(Z_\infty(w)))-\ell_G n|\le 2\eta n\end{array}\right.\end{equation}

Set
$$
r_n= e^{-\ep d(w,\pi_n(Z_{\infty}(w)))}\,.
$$
Therefore, using (\ref{compballs}) with $a=Z_\infty(w)$ and $x=\pi_n (Z_{\infty}(w))$, we get
$$\nu (B_\ep(Z_\infty(w),r_n)) \asymp\nu (\mho_{\GGG}(\pi_n (Z_{\infty}(w)),R))\asymp e^{- d_G(w,\pi_n(Z_{\infty}(w)))}$$
where the right-hand part comes from the fact that $\nu$ is a  \qc measure
of dimension $1/\ep$ for the Green visual metric  and the lemma of the shadow (Lemma \ref{ashadow}).
Hence we deduce  from (\ref{disdri}) that, if $n$ is large enough, then
\begin{equation}\label{firstlimit}
\left|\frac{\log \nu (B_\ep(Z_\infty(w),r_n)) }{\log r_n} -\frac{\ell_G}{\ep\ell} \right|\lesssim \eta
\,.
\end{equation}

Since the measure $\nu$ is doubling (Proposition \ref{cdoub}), $\nu$ is also $\al$-homogeneous
for some $\al >0$,
(cf. \cite[Chap.\,13]{hei})
i.e., there is a constant $C>0$ such that, if $0<r<R<\diam \partial X$ and $a\in\partial X$, then
$$\frac{\nu(B_\ep(a,R))}{\nu(B_\ep(a,r))}\le C\left(\frac{R}{r}\right)^\al\,.$$
 From $$\left|\log \frac{e^{- \ep n\ell} }{r_n} \right|\le 2n\ep\eta$$
it follows that
$$\left|\log \frac{\nu(B_\ep(Z_\infty(w),e^{- \ep n\ell} ))}{\nu(B_\ep(Z_\infty(w),r_n))} \right|
\le 2n\al\ep\eta +O(1)\,.$$
Therefore
$$ \limsup_n \left|\frac{ \log \nu (B_\ep(Z_\infty(w),e^{-\ep n\ell}))}{\log e^{-\ep n\ell}}-
\frac{\log \nu (B_\ep(Z_\infty(w),r_n)) }{\log r_n}\right| \lesssim\eta\,.$$

Since $\eta>0$ is arbitrary, it follows from (\ref{firstlimit}) that
$$\lim_{r\to 0}\frac{ \log \nu (B_\ep(Z_\infty(w),r))}{\log r}=
\lim_{n\to\infty}\frac{ \log \nu (B_\ep(Z_\infty(w),e^{-\ep n\ell}))}{\log e^{-\ep n\ell}}=
\lim_{n\to\infty}\frac{ \log \nu (B_\ep(Z_\infty(w),r_n))}{\log r_n} = \frac{\ell_G}{\ep \ell}\,.$$

In other words,  for $\nu$ almost every $a\in\partial X$,
$$\lim_{r\to 0} \frac{ \log \nu (B_{\ep}(a,r))}{\log r} = \frac{\ell_G}{\ep\ell}\,.$$ \endp

It remains to prove that $\nu$ has dimension $\ell_G/\ep\ell$. This is standard.

{\noindent\bf Hausdorff measures.}
Let $s,t \ge 0$, we set
$$\HHH_s^t(X) \fdefeq \inf\left\{\sum r_{i}^s,\ B_i=B(x_i,r_i),\ X\subset (\cup B_i), r_i\le t\right\}\,, $$
where we consider covers by balls.

The $s$-dimensional measure is then $$\HHH_s(X) \fdefeq \lim_{t\to 0} \HHH_s^t(X)=\sup_{t>0} \HHH_s^t(X)\,.$$

The {\it Hausdorff dimension} $\dim_H X$ of $X$ is the number $s\in[0,\infty]$ such that, for $s'<s$,  $\HHH_{s'}(X)=\infty$ holds
and for all $s'>s$, $\HHH_{s'}(X)=0$.

The {\it Hausdorff dimension} $\dim \nu$ of a measure $\nu$
is the infimum of the Hausdorff dimensions over all sets of full measure.

Replacing covers by balls by covers by any kind of sets in the definition of $\HHH_s^t(X)$
and replacing radii by diameters would not change the value of $\dim \nu$.

For more properties, one can consult \cite{mat}.

\REFPROP{dimH} Let $X$ be a proper metric space and $\nu$ a Borel regular
probability measure on $X$. If, for $\nu$-almost every $x\in X$,
$$\lim_{r\to 0} \frac{\log \nu(B(x,r))}{\log r}=\al$$ then
$\dim\nu =\al$.\ENDPROP

We recall the proof for the convenience of the reader. We will use
the following covering lemma.

\REFLEM{5r} Let $X$ be a proper metric space and $\BBB$ a family of balls
in $X$ with uniformly bounded radii. Then there is a subfamily
$\BBB'\subset\BBB$ of pairwise disjoint balls such that $$\cup_{\BBB} B\subset \cup_{\BBB'} (5B)\,.$$
\ENDLEM

For a proof of the lemma, see \Th 2.1 in \cite{mat}.

{\noindent\sc Proof of Prop.\,\ref{dimH}.} Let $s>\al$, and choose $\eta >0$ small enough so that
$\be:=s-\al-\eta >0$. For $\nu$-almost every $x$, a radius $r_x>0$ exists so that
$$\left|\frac{\log \nu(B(x,r))}{\log r}-\al\right|\le \eta\,,$$
for $r\in (0,r_x]$.

Let us denote by $Y=\{x\in X \; :\; r_x<\infty\}$, which is of full measure.
Let us fix $t\in(0,1)$. For any $x\in Y$, we choose $\rho_x=\min\{r_x,t\}$.
We apply Lemma \ref{5r} to $\{B(x,\rho_x)\}$ and obtain a subfamily $\BBB_{t}$.
It follows that $Y$ is covered by $5\BBB_t$ and
$$\begin{array}{lll} \HHH_s^{5t}(Y) & \le \dis\sum_{\BBB_{t}} (5\rho_x)^s  &\le 5^s t^{\be} \dis\sum_{\BBB_{t}} \rho_x^{\al+\eta}\\& \\
& \lesssim t^{\be} \dis\sum_{\BBB_t} \nu (B(x,\rho_x)) &\lesssim t^{\be} \nu\left(\dis\cup_{\BBB_t} B(x,\rho_x)\right)\\ & \\
& \lesssim t^{\be}&\end{array}$$
which tends to $0$ with $t$. Therefore $\HHH_s(Y)=0$ and so $\dim_H Y \le s$ for all $s >\al$. Whence $\dim\nu\le \al$.

Conversely, let $Y$ be a set of full measure. There is a subset $Z\subset Y$ such that   $\nu(Z)\ge 1/2$ and
such that the convergence of $\log \nu(B(x,r)) / \log r$ to $\al$ is uniform on $Z$  (Egorov theorem). Fix $s<\al$ and let us consider $\eta>0$ small
enough so that $\g=\al-\eta -s >0$. There exists $0<r_0\leq 5$   such that, for any $r\in (0,r_0)$ and any $x\in Z$,
$$\left|\frac{\log \nu(B(x,r))}{\log r}-\al\right|\le \eta\,.$$

Let $\BBB$ be a cover of $Z$ by balls of radius $\rho_x$ smaller than $t\le r_0/5$.
Pick a subfamily $\BBB_{t}=\{B(x,\rho_x)\}$ using Lemma \ref{5r}.
Then $5\BBB_t$ covers $Z$ and
$$1/2\le \dis\sum_{\BBB_{t}} \nu(5B)  \le 5^{\al-\eta} \dis\sum_{\BBB_{t}}
\rho_x^{\al-\eta}\lesssim t^{\g} \dis\sum_{\BBB_{t}} \rho_x^s\,.$$
This proves that $\HHH_s^t(Z)\gtrsim t^{-\g}$
so that $\dim_H Y\ge \dim_H Z\ge \al$. \endp


\section{Harmonic measure of maximal dimension}

This section is devoted to the proof of \Th \ref{main3} and its corollary.

\subsection{The fundamental equality} We assume that $d\in\DDD(\G)$, $\mu$ is a probability measure with exponential moment such that
$d_G\in\DDD(\G)$. Thus there exists $\la >0$ such that $$E\fdefeq\EE\left[e^{\la d(w,Z_1(w))}\right] <\infty\,.$$

The main issue in the proof of \Th \ref{main3} is the following implication which we prove first:

\REFPROP{maxdim}
 Under the hypotheses of \Th \ref{main3}, if $h=\ell v$, then
$\rho$ and $\nu$ are equivalent.
\ENDPROP

Let $R$ be the constant coming from the lemma of the shadow  (Lemma \ref{ashadow}) and write $\mho(x)$ for $\mho(x,R)$.

Let us now define
$$\vp_n=\frac{\rho(\mho(Z_n(w)))}{ \nu(\mho(Z_n(w)))}
\quad \hbox{and} \quad \phi_n=\log \vp_n\,.$$

Since $\mu^n$ is the law of $Z_n$, observe that, if $\be\in (0,1]$, then
$$\EE [\vp_n^\be]=\sum_{\g\in\G}\mu^n(\g) \left(\frac {\rho(\mho(\g(w)))}{ \nu(\mho(\g(w)))}\right)^\be
\mbox{ and } \EE [\phi_n]=\sum_{\g\in\G}\mu^n(\g) \log \left( \frac {\rho(\mho(\g(w)))}{ \nu(\mho(\g(w)))} \right) \, .$$

We start with two lemmata.

\REFLEM{lmaxdim1}
There are finite constants $C_1\ge 1$ and $\be\in(0,1]$ such that, for all $N\ge 1$,
$$\frac{1}{N}\sum_{1\le n\le N} \EE[\vp_n^\be]\le C_1\,.$$
\ENDLEM

When $\mu$ is finitely supported, one can choose $\be=1$ in the lemma.

\Proof
Let $N\ge 1$ and $1\le n\le N$ be chosen. We will first prove that there are some $\kappa$ and $\be$ independent from $N$ and $n$ such that
\begin{equation}\label{p3}
R_{\kappa}\fdefeq \sum_{\g,\ d(w,\g(w))\ge\kappa N} \left( \frac {\rho(\mho(\g(w)))}{ \nu(\mho(\g(w)))}\right)^{\be}\mu^n(\g) \lesssim 1.
\end{equation}

We have already seen that the logarithmic  volume growth rate for the Green metric is $1$.
Then, from the lemma of the shadow (Lemma \ref{ashadow}) applied to both metrics, we get
 \begin{equation}\label{p1}\nu(\mho(\g(w)))\asymp e^{-d_G(w,\g(w))} = F(w,\g(w)) \asymp G(w,\g (w)) =\sum_k \mu^k(\g) \end{equation}
and \begin{equation}\label{p4} \rho(\mho(\g(w)))\asymp  e^{-vd(w,\g(w))}\,.\end{equation}

On the other hand, since $d_G$ is quasi-isometric to $d$, it follows that there is a constant $c>0$ such that
$$\frac {\rho(\mho(\g(w)))}{ \nu(\mho(\g(w)))}\lesssim e^{c d(w,\g(w))}\,.$$

Hence
$$R_{\kappa}\lesssim \sum_{k\ge \kappa N} e^{c\be k}\sum_{k\le d(w,\g(w))<k+1} \mu^n(\g)\,.$$

But $\mu^n$ is the distribution of $Z_n$ so that
$$\sum_{k\le d(w,\g(w))<k+1} \mu^n(\g)\le \PP[d(w,Z_n(w))\ge k]\,.$$
From the exponential Tchebychev inequality, one obtains
\begin{equation}\label{p5}
R_\kappa \lesssim \sum_{k\ge \kappa N} e^{(c\be-\la) k}\EE\left[e^{\la d(w,Z_n(w))}\right]\end{equation}

Now,
$$d(w,Z_n(w))\le \sum_{0\le j<N} d(Z_j(w),Z_{j+1}(w))=\sum_{0\le j<N} d(w,Z_j^{-1}Z_{j+1}(w))$$
since $\G$ acts by isometries. Thus, the independance of the increments of the walk implies
 $$\EE\left[e^{\la d(w,Z_n(w))}\right] \le E^N\,.$$

If we take $\be\fdefeq \min \{\la/2c,1\}$ then (\ref{p5}) becomes
$$R_\kappa\lesssim \sum_{k\ge \kappa N} e^{(-\la/2) k}E^N\lesssim  e^{-(\la/2) \kappa N}E^N\,.$$
The estimate (\ref{p3}) is obtained by choosing $\kappa=2\log E/\la$.

\bigskip

We now prove
that
\begin{equation}\label{p6}
P_N\fdefeq\frac{1}{N}\sum_{1\le n\le N}\sum_{\g\in\G:d(w,\g(w))\le\kappa N} \left( \frac {\rho(\mho(\g(w)))}{ \nu(\mho(\g(w)))}\right)^{\be}\mu^n(\g) \lesssim 1\,.\end{equation}
Both (\ref{p3}) and (\ref{p6}) implies the lemma.

Note that since $\be\le 1$, it follows that $\vp_n^\be\le \max\{1,\vp_n\}\le 1+\vp_n$.

Hence:$$
\begin{array}{ll}
P_N & \lesssim 1 + \dis\frac {1}{ N}\sum_{n=1}^N \dis\sum_{\g\in\G:d(w,\g(w))\le\kappa N}
\dis\frac {\rho(\mho(\g(w)))}{\nu(\mho(\g(w)))} \mu^n(\g)\\ & \\
& \lesssim 1+  \dis\frac {1}{ N} \dis\sum_{\g\in\G:d(w,\g(w))\le\kappa N}
\dis\frac {\dis\sum_{n=1}^N\mu^n(\g) }{\nu(\mho(\g(w)))} \rho(\mho(\g(w))).\end{array}
$$
But (\ref{p1}) implies that $$\frac {\dis\sum_{n=1}^N\mu^n(\g) }{\nu(\mho(\g(w)))}\lesssim 1$$
so that
\begin{equation}\label{p2} P_N\lesssim 1 + \dis\frac {1}{ N} \dis\sum_{d(w,x)\le\kappa N}
\rho(\mho(x))\,.\end{equation}

Since $\rho(\mho(x))\asymp e^{-vd(w,x)}$ by (\ref{p4}) and since there are
approximately $e^{vk}$ elements in the $d$-ball of radius $k$ (Theorem \ref{reg}),
we have
$$\sum_{d(w,x)\le \kappa N} \rho(\mho(x))\asymp  \sum_{1\le n \le \kappa N}  e^{vn}e^{-vn}\,,$$ and
$$\sum_{d(w,x)\le \kappa N} \rho(\mho(x))\lesssim N\,.$$

Therefore, the estimate (\ref{p6}) follows from (\ref{p2}).
\endp

\REFLEM{lmaxdim2}
There is a finite constant $C_2 \ge 0$  such that the sequence $(\EE(\phi_n)+C_2)_{n\ge 1}$
is subadditive and $(1/n)\phi_n$ tends to $h-\ell v$ a.s. and in expectation.
\ENDLEM

\Proof
By the lemma of the shadow (Lemma \ref{ashadow}),
$$\frac{1}{n}\phi_n= \frac{1}{n}d_G(w,Z_n(w))- \frac{1}{n}v d(w,Z_n(w)) +O(1/n)$$
so, from Kingman ergodic \th it follows that $(1/n)\phi_n$ converges almost
surely and in expectation towards $$\ell_G-\ell v= h-\ell v\,,$$
since  $h=\ell_G$, see \cite{bhm}.

Let $m,n\ge 1$. It also follows from the lemma of the shadow and the triangle inequality for $d_G$ that
$$
\EE[\phi_{m+n}]- (\EE[\phi_m]+\EE[\phi_n]) \le
v \EE[d(w,Z_m(w))+d(Z_m(w),Z_{m+n}(w))-d(w,Z_{m+n}(w))]+O(1)\,.$$
So Proposition \ref{control} implies the existence of some constant  $C_2$ such that
$$
\EE[\phi_{m+n}]- (\EE[\phi_m]+\EE[\phi_n]) \le C_2\,.$$
This gives the desired subadditivity.
\endp

{\noindent \sc Proof of Proposition \ref{maxdim}.}
We shall prove that if $\rho$ and $\nu$
are not equivalent, then $h<\ell v$.

Assuming that $\rho$ and $\nu$ are not equivalent, the ergodicity of both measures implies that
$\vp_n$ tends to $0$ $\PP$-a.s.

Choose $\eta\in (0,e^{-1}]$.

By Egorov theorem, there exist two measurable sets $A$ and $B=A^c$ such that
$\PP [A]\le \eta$ and $(\vp_n|_B)_n$ converges uniformly to $0$.

For any $n\ge 1$,
$$\EE[\phi_n]= \int_A \phi_n d\PP + \int_B\phi_n d\PP\,.$$

Since $(\vp_n|_B)_n$ uniformly converges to $0$, there exists $n_0$ such that for $n$ larger than $n_0$,
$\phi_n|_B\le \log \eta$ and therefore
$$\int_B\phi_n d\PP\le \PP [B]  \log \eta \le (1-\eta)\log \eta\,.$$

Choose $\be$ and $C_1$ as in Lemma \ref{lmaxdim1}.
Jensen inequality yields
$$\int_A \phi_n d\PP \le \frac{\PP[A]}{\be}\log \int_A\vp_n^\be \frac{d\PP}{\PP[A]}\le \frac{\eta}{\be}\log(1/\eta) +\frac{\eta}{\be} \log^+\EE[\vp_n^\be]\,,$$
where we have used $\eta\le 1/e$.

But Lemma \ref{lmaxdim1} implies that $\liminf \EE[\vp_n^\be] < 2C_1$. So that there exists $p\ge n_0$
such that $\EE[\vp_p^\be]\le 2C_1$.

Hence,
$$\EE[\phi_p] \le   (1-\eta) \log \eta + \frac{\eta}{\be}\log(1/\eta) + \frac{\eta}{\be}\log (2C_1)\,.$$

When $\eta$ tends to $0$, the right-hand side tends to $-\infty$. Therefore, if we fix $\eta$ small enough,
 there exists $p$ such that

$$
\EE[\phi_p] + C_2\le -1
\,,
$$
where $C_2$ is the constant appearing in Lemma \ref{lmaxdim2}.

Lemma \ref{lmaxdim2} now implies that
$$\frac{1}{k} (\EE[\phi_{kp}]+ C_2)\le \EE[\phi_p]+C_2\le -1$$
for $k\ge 1$.
As $(1/pk) \EE[\phi_{pk}]$ tends to $(h-\ell v)$, letting $k$ go to infinity, one obtains
$$(h-\ell v)\le \frac{-1}{p}<0\,.$$\endp

{\noindent\bf Remark.}
 In view of the proof of Proposition \ref{maxdim}, one might
wonder whether it is always true that a
doubling measure of maximal dimension in an Ahlfors-regular space,
as $\nu$ is, has to be equivalent to the Hausdorff measure of the
same dimension. This property turns out to be false in general.
We are grateful to P.\,Mattila for pointing out to us its invalidity
and to Y.\,Heurteaux for providing an explicit example of a doubling measure of dimension
$1$ in the unit interval $[0,1]$ which is singular to the Lebesgue measure.
We briefly describe his construction.

Let us consider a sequence of integers $(T_n)_n$ tending to infinity and satisfying
\ben
\item   $T_{2n}-T_{2n-1}$ is equivalent to $T_{2n}$;
\item $T_{2n+1}-T_{2n}$ is negligible in front of $T_{2n}$;
\item $T_{2n+1}-T_{2n}$ tends to infinity.\een

We then fix a weight $p\in (0,1/2)$, 
define the sequence $(p_k)_k$ as follows:
\bit
\item If  $T_{2n-1}< k\le T_{2n}$, then $p_k=1/2$;
\item If  $T_{2n}< k\le T_{2n+1}$, then $p_k=p$. \eit

This means that the value of  $p_k$ is either $p$ or $ 1/2$, and that, asymptotically, the mean of
$( p_k)_{1\le k\le n}$ tends towards $1/2$: $\lim_N (1/N) \sum_{k=1}^N p_k=1/2$.

Let us code the dyadic intervals of $[0,1]$ of the $n$th  generation as $I_{a_1\ldots a_n}$, with  $a_1,\ldots \,,a_n =0$ or $1$.
We  define the measure $m$ by setting

$$\frac{m(I_{a_1\cdots a_{n+1}})}{m(I_{a_1\cdots a_n})}=p_n\quad \hbox{if} \quad a_{n+1}=a_n$$

and

$$\frac{m(I_{a_1\cdots a_{n+1}})}{m(I_{a_1\cdots a_n})}=1-p_n\quad \hbox{if} \quad a_{n+1}\not= a_n\,.$$

The measure $m$  we have just defined has dimension $1$ (any set of dimension less than $1$ is $m$-negligible),
is doubling (see \cite{heur} for similar constructions),
but it is singular with respect to the Lebesgue measure since
$2^n m(I_{a_1\ldots a_n})$ tends to $0$ a.e.

\subsection{Equivalent measures}

We let $\G$ be a non-elementary hyperbolic group, $(X,d)\in\DDD(\G)$, and $\mu$ a probability measure
on $\G$ so that $d_G\in\DDD(\G)$. This section is devoted to proving
\REFPROP{equivmeasure} If $\rho$ and $\nu$ are equivalent then their density is almost surely bounded
i.e., there is a constant $C\ge 1$ such that
for any Borel set $A\subset\partial X$,
$$\frac{1}{C}\nu(A)\le \rho(A)\le C\nu (A)\,.$$
\ENDPROP

We will work with the space $\partial^2 X$ of distinct points $(a,b)\in\partial X\times\partial X$, $a\ne b$,
which is reminiscent to the geodesic flow of a negatively curved manifold. The group $\G$ acts on $\partial^2 X$ by
the diagonal action $\g\cdot(a,b)=(\g(a),\g(b))$, $\g\in\G$.

We define the following two $\si$-finite measures on $\partial^2 X$\,:
$$d\tilde{\rho}(a,b)= \frac{d\rho(a)\otimes d\rho(b)}{\exp 2v(a|b)} \quad\hbox{and} \quad d\tilde{\nu}(a,b)=\frac{d\nu(a)\otimes d\nu(b)}{\exp 2(a|b)^G},$$
where we define $$(a|b)^G\fdefeq \liminf_{(a_n),(b_n)\to a,b} (a_n|b_n)^G\,.$$
We recall that since $\nu$ is a conformal measure, $\tilde{\nu}$ is invariant, and
it is furthermore ergodic \cite[Thm 3.3]{kai94}.
On the other hand, $\rho$ being just a \qc measure,
it follows that $\tilde{\rho}$ is just quasi-invariant, cf. \cite{co}. This implies the
existence of a constant $C\ge 1$ such that, for any Borel set $A\subset\partial^2 X$,
$$\frac{1}{C}\tilde{\rho}(A)\le \tilde{\rho}(\g(A))\le C\tilde{\rho} (A)\,.$$

\noindent{\sc Proof of Proposition \ref{equivmeasure}.}
By assumption, there is a positive $\nu$-integrable function $J$ such that
$d\rho =J d\nu$.
Therefore,
$d\tilde{\rho} =\tilde{J} d\tilde{\nu}$ holds with $$\tilde{J}(a,b)=
J(a)J(b)\frac{\exp 2(a|b)^G}{\exp 2v(a|b)}\, .$$

We shall first prove that $\tilde{J}$ is essentially constant (and non-zero).
There is a constant $C>1$ such that the set $$A\fdefeq\{ (1/C)\le \tilde{J}\le C\}$$
has positive $\tilde{\nu}$-measure.
Since $\tilde{\nu}$ is ergodic, for $\tilde{\nu}$-almost every $(a,b)\in\partial^2 X$,
there exists $\g\in\G$ such that $\g(a,b)\in A$.
It follows from the invariance of $\tilde{\nu}$ and the quasi-invariance of $\tilde{\rho}$ that
$$\tilde{J}(a,b)\asymp \tilde{J}(\g(a),\g(b))\,.$$
This proves the claim.

Therefore, for $\tilde{\nu}$-almost every $(a,b)$,
$$J(a)J(b)\asymp \frac{\exp 2v(a|b)}{\exp 2(a|b)^G}\,.$$
Let us assume that $\log J$ is unbounded in a \nbhd $U$ of a point $a\in\partial X$.
We may find a point $b\in\partial X$ with $J(b)$ finite and non-zero, and far enough from $U$
so that $$\frac{\exp 2v(c|b)}{\exp 2(c|b)^G}\asymp 1$$
for any $c\in U$. This proves that $\log J$ had to be bounded in $U$\,: a contradiction. \endp

\subsection{Geometric characterisation of the fundamental inequality}
We may now turn to the proof of \Th \ref{main3}.

{\noindent\sc Proof of \Th \ref{main3}.} We first prove that (i), (ii) and (iii) are equivalent.
Then we prove that  (iii) implies (iv),  (iv) implies (v) which implies (iii).
\bit
\item From Proposition \ref{maxdim}, we deduce that (i) implies (ii). Proposition \ref{equivmeasure} says that (ii) implies (iii).
Furthermore, if $\nu$ and $\rho$ are equivalent,
then they have the same Hausdorff dimension. So, from Corollary \ref{cmain2}
and \Th \ref{reg2}, we get
that
$$
\frac{h}{\ell \ep} = \dim \nu = \dim \rho = \frac{v}{\ep} \, ,
$$
and thus $h= \ell v$.

\item To prove  that (iii) implies (iv), we apply the lemma of the shadow (Lemma \ref{ashadow}): it follows that,
for any $\g\in\G$, $$e^{-v d(w,\g(w))}\asymp\rho(\mho(\g(w)))\asymp \nu(\mho(\g(w)))\asymp e^{-d_G(w,\g(w))}$$
whence the existence of a constant $C$ such that
$$|vd(w,\g(w))-d_G(w,\g(w))|\le C\,.$$ Since $\G$ acts  transitively by isometries for both metrics, it follows
that $(X,vd)$  and $(X,d_G)$ are $(1,C)$-quasi-isometric.

\item Assuming (iv), it follows that Busemann functions coincide up to the  multiplicative factor $v$.
Therefore, the Radon-Nikodym derivative of $\g^*\nu$ with respect to
$\nu$ at a point $a\in\partial X$ is  proportional to $\exp (-v\be_a(w,\g^{-1}(w)))$ a.e.
Therefore, $\nu$ is a \qc measure for $(\partial X, d_\ep)$. This is (v).

\item For the last implication, (v) implies (iii),
one can use the uniqueness statement in \Th \ref{reg}
to get that $\rho$ and $\nu$ are equivalent and have bounded density. This proves (iii).\endp \eit

\subsection{Simultaneous random walks} We now turn to the proof of Corollary \ref{cmain3}.

{\noindent\sc Proof of Corollary \ref{cmain3}.} Let us consider the Green metric $d_G$ associated with
$\mu$ and denote by $\hat{\ell}$ the drift of $(\hat{Z}_n)$ in the metric space $(\G,d_G)$.
\Th \ref{main1} implies that $d_G\in\DDD(\G)$.

Assumption (i) translates into $\hat{h}=\hat{\ell}$. Since $v_G=1$, this means that $\hat{\nu}$
has maximal dimension in the boundary of $(\G,d_G)$ endowed with a visual metric.
Therefore \Th \ref{main3} implies the equivalence between (i) and (iii).

Exchanging the roles of $\mu$ and $\hat\mu$ gives the equivalence between (ii) and (iii).

If $\hat{d}_G$ denotes the Green metric for $\hat{\mu}$,
then  (iv)  means that $d_G$ and $\hat{d}_G$ are $(1,C)$-quasi-isometric,
which is equivalent to (iii) by \Th \ref{main3}.
\endp

\subsection{Fuchsian groups with cusps} 
 
The next proposition is the key to the proof of Theorem \ref{GuivLeJan}. 
Let us first introduce its setting.

Let $X$ be a proper quasiruled hyperbolic space and let $\Gamma$ be a hyperbolic subgroup of isometries  
that acts properly discontinuously on $X$. Consider a symmetric probability measure  $\mu$ 
on $\Gamma$ with finite support and whose support generates the group $\Gamma$. Let $\nu$ be the 
corresponding harmonic measure on $\partial\Gamma$, the visual boundary of $\Gamma$. 

Let $\Gamma(w)$ be the orbit of some point $w\in X$. 
As for Theorem \ref{kai}, Theorem 7.3 in \cite{kai} and \S 7.4
therein imply that 
the sequence $Z_n(w)$ almost surely 
converges to some point $Z_\infty(w)$ in $\partial X$, the visual boundary of $X$. Let $\nu_1$ be the 
law of $Z_\infty(w)$. 

Although the two spaces $\partial X$ and $\partial \Gamma$ might be topologically different, the two 
measured spaces $(\partial X,\nu_1)$ and $(\partial \Gamma,\nu)$ are isomorphic as $\Gamma$-spaces  
i.e., there exists a measured spaces isomorphism $\Phi$ from $(\partial \Gamma,\nu)$ to 
$(\partial X,\nu_1)$ that conjugates the action of $\Gamma$ on both spaces. 
Indeed both spaces are models for the Poisson boundary of the random walk. 
This is proved in \cite{kai} Theorem 7.7 and Remark 3 following it for $(\partial X,\nu_1)$ and 
it is a general fact for the Martin boundary $(\partial \Gamma,\nu)$.

\REFPROP{equivmeasure2} Let $X$ be a proper quasiruled hyperbolic space endowed with a geometric group action. 
Let $\rho$ be the corresponding Patterson-Sullivan measure. Let $\Gamma$ be a hyperbolic subgroup of isometries  
that acts properly discontinuously on $X$ and $\Gamma(w)$ be an orbit of $\Gamma$ in $X$. 
Let $\mu$ be a symmetric probability measure on $\Gamma$ 
with finite support and whose support generates the group $\Gamma$. Let $\nu_1$ be the limit law 
of the trajectories of the random walk on $\partial X$. 
If $\rho$ and $\nu_1$ are equivalent then $\Gamma$ and $\Gamma(w)$ are quasi-isometric. 
\ENDPROP

\Proof 
One checks as in Proposition \ref{equivmeasure} that, once $\rho$ and $\nu_1$ are equivalent, 
then their density is almost surely bounded. 

We recall the 
following change of variables formula: 
$$ \frac{d \g^*\nu}{d\nu}(a)=K_a(\g^{-1})\,,$$
for $\nu$ almost any point $a\in\partial\Gamma$ and where $K_a$ is the Martin kernel. 
Because of the isomorphism $\Phi$, we also have 
$$ \frac{d \g^*\nu_1}{d\nu_1}(\xi)=K_{\Phi^{-1}(\xi)}(\g^{-1})\,,$$
for $\nu_1$ almost any point $\xi\in\partial X$.  

On the other hand, $\rho$ being a quasiconformal measure, it satisfies 
$$\frac{ d \g^*\rho}{d\rho}(\xi)\asymp e^{v\be_\xi(w,\g^{-1}(w))} \,,$$ 
where $\be_\xi$ is the Busemann function in $X$. 

Since the density of $\nu_1$ with respect to $\rho$ is bounded and bounded away from $0$, 
we therefore have 
\begin{equation} \label{equiv}
K_{\Phi^{-1}(\xi)}(\g^{-1})\asymp e^{v\be_\xi(w,\g^{-1}(w))}\,,
\end{equation} 
for $\rho$ almost any $\xi$. 

We now use Lemma \ref{aa4}. First observe that $\sup_{\xi\in\partial X} \be_\xi(x,y)$ can be replaced 
by an essential sup with respect to $\rho$ since $\rho$, being quasiconformal, charges any non 
empty ball and since $\xi\rightarrow\be_\xi(x,y)$ is locally almost constant. 
So we get from Lemma \ref{aa4} that 
$$ \vert  d(x,y)- v\, ess\!\sup_{\xi\in\partial X} \be_\xi(x,y) \vert $$ 
is bounded. By a similar argument, applying Lemma \ref{aa4} to the Green metric on $\Gamma$, 
we deduce that 
$$ \vert d_G(e,\gamma^{-1}) - ess\!\sup_{\xi\in\partial \Gamma} \log K_\xi(\gamma^{-1}) \vert $$ 
is bounded. The essential sup is taken with respect to $\nu$. 

But (\ref{equiv}) implies that 
$$\vert ess\sup_{a\in\partial \Gamma} \log K_a(\gamma^{-1}) - 
v\, ess\!\sup_{\xi\in\partial X} \be_\xi(w,\gamma^{-1}(w)) \vert $$ 
is bounded and therefore 
$$\sup_{\gamma\in\Gamma} \vert d_G(e,\gamma^{-1})-v\, d(w,\gamma^{-1}(w)) \vert<\infty\,.$$ 
We conclude that $\Gamma$ and $\Gamma(w)$ are indeed quasi-isometric.  

\endp

{\noindent\sc Proof of Theorem  \ref{GuivLeJan}.} 
We proceed by contradiction and assume that $\nu_1$ is equivalent to 
the Lebesgue measure $\lambda$ on $\SS^1$. 

First note that we can restrict our attention to the subgroup generated by the support of $\mu$. 
If this subgroup turned out to have infinite covolume then its boundary would be a strict subset of 
$\SS^1$ and $\nu_1$ would certainly not be equivalent to the Lebesgue measure. Therefore 
we may, and will, assume that the support of $\mu$ generates $\Gamma$ 
and that $\Gamma$ has finite covolume and is finitely generated. 

We know from Selberg's lemma that $G$ contains a torsion-free finite subgroup $\G_S$
of finite index so that $\HH^2/\G_S$ is a compact Riemann surface with finitely many
punctures. Therefore $\G_S$ is isomorphic to a free group so that $\G$ is  
hyperbolic and its 
boundary is a perfect, totally disconnected, compact set (a Cantor set). 

Let $\Gamma(w)$ be an orbit of $\Gamma$ in $\HH^2$. 
By the finite covolume assumption, the limit set of $\Gamma(w)$ is homeomorphic to the circle $\SS^1$. 
But it follows from Proposition \ref{equivmeasure2}  
that $\Gamma$ is quasi-isometric to $\Gamma(w)$. As a consequence 
the limit set of $\Gamma(w)$ is  also homeomorphic to  the boundary of $\Gamma$. 
As $\SS^1$ is not a Cantor set, we get the contradiction we were looking for.   
\endp


\section{Discretisation of Brownian motion}\label{sec:bm}
We let $M$ be the universal covering of a
Riemannian manifold $N$ of  pinched negative curvature
and finite volume
with deck transformation group $\G$ i.e.,
$M/\G=N$. We let $d$ denote the distance defined
by the Riemannian structure on $M$.
Note that when $N$ is compact,  $\G$ acts geometrically on $M$, and since it has negative curvature,
it follows that $\G$ is hyperbolic and that $M$ is quasi-isometric to $\G$
by \v{S}varc-Milnor's lemma (Lemma \ref{svmi}).

We  consider the diffusion process $(\xi_t)$ generated by  the Laplace-Beltrami operator $\Delta$  on
$M$.
That is, we let $p_t$ be the fundamental solution of the heat
equation $\partial_t=\Delta$. Then there is a probability measure
$\PP^y$ on the family $\Xi^y$ of continuous curves $\xi:\R_+\to M$ with $\xi_0=y$ such that,
for any Borel  sets $A_1$, $A_2$, ..., $A_n$, and any times
$t_1<t_2<\ldots <t_n$,
\begin{eqnarray*}
&&\PP^y(\xi_{t_1}\in A_1,\ \ldots, \xi_{t_n}\in A_n)\\
&&=
\int_{A_1}\int_{A_2}\ldots\int_{A_n}p_{t_1}(y,x_1)p_{t_2-t_1}(x_1,x_2)\ldots
p_{t_n-t_{n-1}}(x_{n-1},x_n)dx_1\ldots dx_n\,.\end{eqnarray*}
If $\mu$ is a positive measure on $M$, we write $\PP^\mu=\int_M\PP^y\mu(dy)$,
and this defines
a measure on the set of Brownian paths $\Xi$.

As for random walks, the following limit exists almost surely and in $L^1$ and we call it
the drift of the Brownian motion:
$$\ell_M\fdefeq\lim \frac{d(\xi_0,\xi_t)}{t}\,;$$
it is also known that $\ell_M>0$ and that $(\xi_t)$ almost surely converges to a point $\xi_\infty$ in $\partial M$
\cite{pra,pin}.
The distribution of $\xi_\infty$ is the harmonic measure.
Furthermore, V.\,Kaimanovich has defined an asymptotic entropy $h_M$ which shares the
same properties as for random walks \cite{kai86}: for any $y\in M$,
$$h_M\fdefeq\lim \frac{-1}{t}\int p_t(y,x)\log p_t(y,x)dx \, .$$

He also proved that the fundamental inequality $h_M\le\ell_M v$ remains valid in this setting,
where $v$ denotes the logarithmic volume growth rate of $M$.

\subsection{The discretised motion}
W.\,Ballmann and F.\,Ledrappier have refined a method
of T.\,Lyons and D.\,Sullivan \cite{ls}, further studied by A.\,Ancona \cite{aa},
V.\,Kaimanovich \cite{kai92}, and by A.\,Karlsson and F.\,Ledrappier \cite{kl07} which replaces the Brownian motion by a random walk
on $\G$ \cite{bl}. The construction goes as follows in our specific case.

Let $\pi:M\to N$ be the universal covering and let us fix a base
point $w\in M$. Fix $\ep>0$ smaller than the injectivity radius
of $N$ at $\pi(w)$, and consider $\overline{V}=B(\pi(w),\ep)$ in $N$;
for $D$ large enough, the set  $\overline{F}=\{G_{\overline{V}}(\pi(w),\cdot) \ge D\}$
is compact in $\overline{V}$, where $G_{\overline{V}}$ denotes the Green function
of the Brownian motion killed outside $\overline{V}$.
There exists a so-called Harnack constant $C<\infty$ such that, for any positive harmonic function
$h$ on $\overline{V}$ and any points $a,b\in \overline{F}$, $h(a)/h(b)\le C$ holds.

Let $V=\pi^{-1}(\overline{V})$, $F=\pi^{-1}(\overline{F})$,  $V^x=B(x,\ep)$ and $F_x=F\cap V^x$
for $x\in X\fdefeq\G(w)$. If $y\in F_x$, we set $\chi(y)=x$.
(Note that $\chi$ is well defined thanks to the choice of $\ep$.)

Let $\xi_t$ be a sample path of the Brownian motion. We define inductively the following
Markov stopping times $(R_n)_{n\ge 1}$ and $(S_n)_{n\ge 0}$ as follows.

Set $S_0=0$ if $\xi_0\notin X$, and $S_0=\min\{t\ge 0,\ \xi_t\notin V^{\xi_0}\}$.
Then, for $n\ge 1$, let
$$\left\{\begin{array}{l} R_n= \min\{t\ge S_{n-1},\ \xi_t\in F\}\\
S_n= \min\{t\ge R_n,\ \xi_t\notin V^{X_n}\}\end{array}\right.$$
with $X_n=\chi(\xi_{R_n})$.

Let us also define recursively for $k\ge 0$ on $\Xi\times [0,1]^{\N}$,
$$\left\{\begin{array}{l}N_0(\xi,\al)=0\\
\\
N_k(\xi,\al)=\min\{n> N_{k-1}(\xi,\al),\ \al_n <\kappa_n(\xi)\}\end{array}\right.$$
where $$\kappa_n(\xi)=\frac{1}{C}\frac{d\ep_{X_n}^{V}}{d\ep_{\xi_{R_n}}^V}(\xi_{S_n})\,,$$
and, for $z\in F$, $\ep_z^V$ denotes the distribution of $\xi_{S_1}$ for sample paths $\xi_t$ starting at $z$. We also set
$$T_k=S_{N_k}\,.$$

For $y$ in $M$, we let $\tP^y$ denote the product measure of $\PP^y\times \la^{\N}$, where $\la$
is the Lebesgue measure on $[0,1]$. We then define on $X$, the law
$$\mu_y(x)=\tP^y[X_{N_1}=x]\,.$$

The following properties are known to hold \cite{ls,kai92,bl,kl07}.

\REFTHM{propdis} Let us define $\mu(\g)=\mu_w(\g(w))$, and  $Z_k(w)=X_{N_k}$ with $Z_0(w)=w$.
\ben
\item[(i)]  The random sequence $(Z_n(w))$ is the random walk generated by $\mu$: for any $x_1=\g_1(w),\ldots,x_n=\g_n(w)\in X$,
$$\tP^w(Z_1=x_1,\ldots,Z_n=x_n)=\mu(\g_1)\mu(\g_1^{-1}\g_2)\ldots \mu(\g_{n-1}^{-1}\g_n)\,.$$
\item[(ii)] The measure $\mu$ is symmetric with full support but has a finite  first
moment with respect to $d$.
\item[(iii)] The Green function $G_\mu$ of the random walk is proportional
to the Green function $G_M$ of $M$.
\item[(iv)] There exists a positive constant $T$ such that the following limit exists almost surely and in $L^1$:
$$\lim\frac{S_{N_k}}{k}=T\,.$$
\item[(v)] Almost surely and in $L^1$, $$\lim \frac{d(\xi_{kT},Z_k(w))}{k}=0\,.$$
\item[(vi)] The harmonic measures for the Brownian motion and the random walk coincide.\een
\ENDTHM

We are able to prove the following:

\REFTHM{hyp+dim} Under the notation and assumptions from above,
let $d_G$ denote the Green metric associated with  $\mu$.  If $N$ is compact, then
$d_G\in\DDD(\G)$ and
$$\dim \nu= \frac{h_M}{\ep\ell_M}$$
where $h_M$ and $\ell_M$ denote the entropy and the drift of the Brownian motion
respectively.\ENDTHM

\Proof
The acronyms (ED) and (QR) below refer to Proposition \ref{d_GdansD(G)}.
Since $M$ has pinched negative curvature,
it follows that $G_M(x,y)\lesssim e^{-c d(x,y)}$ holds for some constant $c>0$,
see \cite[(2.4) p.\,434]{as}.
By part (iii) of \Th \ref{propdis}, $G_\mu$ and $G_M$ are proportional.
Therefore $G_\mu$ also satisfies $G_\mu(x,y)\lesssim e^{-c d(x,y)}$
and (ED) is proved.
Furthermore, A.\,Ancona's \Th \ref{anc} also holds for the Brownian motion,
see \cite{aa}, showing that (QR) holds as well.
Both these properties imply that $(X,d_G)\in\DDD(\G)$ by Proposition \ref{d_GdansD(G)}.

The identity $h_\mu=h_M\cdot T$ was proved by V.\,Kaimanovich \cite{kai86,kai92}.
Furthermore, from \Th \ref{propdis} (4), it follows that almost surely,
$$\ell_\mu =\lim \frac{d(w,Z_k(w))}{k}=\lim \frac{d(w,\xi_{kT})}{k}=\ell_M\cdot T\,.$$
Thus, Corollary \ref{cmain2} implies that
$$\dim\nu = \frac{h_M}{\ep \ell_M}\,.$$\endp

The computation of the drift can also be found in \cite{kl07}.

\subsection{Exponential moment for the discretised motion}
In \cite{aa}, A.\,Ancona  wrote in a remark that
the random walk defined above has a finite exponential moment  when $N$ is compact.
Since this fact is crucial to us, we provide  here a detailed proof.
This will enable us to apply \Th \ref{main3} and conclude the proof
of \Th \ref{main4}.

\REFTHM{finitexp}  If $N$ is compact, then the  random walk  $(Z_n)$ defined in \Th \ref{propdis} has a finite exponential
moment.\ENDTHM

The proof requires intermediate estimates on the Brownian
motion.
The main step is an estimate on the position of $\xi_{S_1}$:

\REFPROP{xi_S1>r} There are positive  constants $C_1$ and $c_1$ such that, for any $r\ge 1$,
$$\sup_{y\in M} \PP^y[d(\xi_0,\xi_{S_1})\ge r] \le C_1 e^{-c_1 r}\,.$$\ENDPROP

Proposition \ref{xi_S1>r} follows from the following lemma.

\REFLEM{supxi} We write $\xi_t^*=\sup_{0\le s\le t}d(\xi_0,\xi_s)$.
There are constants $m>0$, $c_2>0$ and $C_2>0$ such that
$$\sup_{y\in M}\PP^y\left[\xi_t^*\ge m t\right]\le C_2 e^{-c_2 t}\,.$$
\ENDLEM

\Proof
We first prove that all the exponential moments of $\xi_1^*$ are finite.
Our proof relies on the following upper Gaussian estimate valid as soon as the curvature
is bounded
(see e.g. \cite[\S\,6]{pra} for a proof): for any $y\in M$ and any $t\ge 2$,
$$\PP^y[\xi_1^*\ge t]\le \exp\left(-c t^2\right)\,,$$
for some constant $c$ that does not depend on $y$ nor on $t$.

Hence, if $\la >0$ then
$$\begin{array}{ll} \EE^y\left[e^{\la\xi_1^*}\right]
& = 1+ \dis\int_{u>0}e^u\PP^y[\xi_1^*\ge (u/\la)]du\\&\\
& \le 1+ \dis\int_{u>0} e^{u-c\frac {u^2}{\la^2}}du <\infty \, . \end{array}$$

\bigskip

Let $y\in M$ and $m>0$.
It follows from the exponential Tchebychev inequality that
$$
\PP^y[\xi_t^*\ge m t] \le e^{-\la mt}\EE^y\left[e^{\la\xi_t^*}\right]\,.$$
 We remark that, for $n\ge 1$ and $t\in (n-1,n]$,
$$\xi_t^*\le \sum_{0\le k< n}\sup_{k\le s\le k+1}d(\xi_k,\xi_s)\,.$$
It follows from the Markov property that, for all $y\in M$,
$$\EE^y\left[e^{\la\xi_t^*}\right]\le\left(\sup_{z\in M}\EE^z\left[e^{\la\xi_1^*}\right] \right)^n\,.$$
Therefore
$$
\PP^y[\xi_t^*\ge m t] \lesssim e^{-\la mt}\left(\sup_{z\in M}\EE^z\left[e^{\la\xi_1^*}\right] \right)^t \, .$$
So, if $m$ is chosen large enough, we will find $c_2>0$ so that
$$\PP^y[\xi_t^*\ge m t] \lesssim e^{-c_2 t}\,.$$

\endp

{\noindent\sc Proof of Proposition \ref{xi_S1>r}.}
The compactness of $N$ easily implies the following upper bound on the first hitting time
$S_1$ using the orthogonal decomposition of $L^2(N)$  (see \cite[(5.2)]{pin}):
there are positive constants $C_3$ and
$c_3$ such that, for any $y\in M$,
\begin{equation}\label{S_1>n}\PP^y[ S_1\ge k] \le C_3e^{-c_3 k}\end{equation}

Let us consider $\kappa >0$ that will be fixed later.
$$\PP^y[d(y,\xi_{S_1})\ge r] \le\PP^y[d(y,\xi_{S_1})\ge r;\,S_1\le\kappa]
+\PP^y[d(y,\xi_{S_1})\ge r;\,S_1\ge \kappa]\,.$$

From (\ref{S_1>n}), it follows that
$$\PP^y[d(y,\xi_{S_1})\ge r] \lesssim \PP^y[\xi_{\kappa}^*\ge r] +  e^{-c_3\kappa}\,.$$

Choosing $\kappa= r/m$, Lemma  \ref{supxi} implies that
$$\PP^y[d(y,\xi_{S_1})\ge r] \lesssim  e^{-\frac{c_2}{m} r}+ e^{-\frac{c_3}{m} r}\,,$$
and the proposition follows.\endp

{\noindent\sc Proof of \Th \ref{finitexp}.} Let $r\ge 1$ and $k\ge 1$, and $\la >0$ that
will be fixed later.

The exponential Tchebychev inequality yields
$$\PP^y[d(\xi_0,\xi_{S_k})\ge r] \le e^{-\la r}\EE^y\left[ e^{\la d(\xi_0,\xi_{S_k})}\right]\,.$$
But $$d(\xi_0,\xi_{S_k})\le \sum_{0\le j< k}d(\xi_{S_j},\xi_{S_{j+1}})$$
so the strong Markov property implies that
\begin{equation}\label{me3}
\PP^y[d(\xi_0,\xi_{S_k})\ge r] \le e^{-\la r}\left(\sup_{z\in M} \EE^z\left[e^{\la d(\xi_0,\xi_{S_1})}\right]\right)^k \, .\end{equation}

Using Proposition \ref{xi_S1>r} and its notation, we get that for any $z\in M$, 

$$\begin{array}{ll} \EE^z\left[e^{\la d(\xi_0,\xi_{S_1})}\right] &
=1+ \dis\int_{u>0}e^u\PP^y[d(\xi_0,\xi_{S_1})\ge (u/\la)]du\\&\\
& \le 1 + C_1 \dis\int_{u>0} e^u e^{-c_1 u/\la}du \, . \end{array}$$
We choose $\la < c_1$; there exists a positive constant $C_4$ such that
$$\sup_{z\in M}\EE^z\left[e^{\la d(\xi_0,\xi_{S_1})}\right]\le \frac{1+C_4\la }{1-(\la/c_1)}\,.$$
Plugging this last inequality in (\ref{me3}) yields
\begin{equation}\label{me4} \PP^y[d(\xi_0,\xi_{S_k})\ge r] \lesssim e^{-\la r +kc_4\la}\end{equation}
for some constant $c_4>0$.

\bigskip

We note that, for any $x\in X$, any $z\in F_x$ and $u\in\partial V^x$,
$$\frac{d\ep_x^V}{d\ep_z^V}(u)\ge (1/C)$$ where $C$ is the Harnack constant.
Observe that this estimate is uniform with respect to $u\in\partial V^x$  and $z\in F_x$.
Therefore,

\begin{eqnarray}
\tP^y[T_1 \ge S_k\vert\,\xi] &  = &\tP^y\left[ \dis\cap_{n=1}^{k-1}\{ \kappa_n(\xi)<\al_n\}\vert\, \xi\right]  \nonumber \\
&\leq & \tP^y\left[ \dis\cap_{n=1}^{k-1} \{(1/C^2)<\al_n\}\vert\, \xi\right] =  \tP^y\left[ \dis\cap_{n=1}^{k-1} \{(1/C^2)<\al_n\}\right]\nonumber \\
 & = & \dis\prod_{n=1}^{k-1} \tP^y[ (1/C^2)<\al_n] \lesssim  (1-(1/C^2))^k\,. \label{me5}
 \end{eqnarray}

 In (\ref{me5}), we used the notation $\tP^y[.\vert\,\xi]$ to denote
the conditional probability given the Brownian path $\xi$.
Note that $S_k$, being a function of $\xi$,
does not depend on the sequence $\alpha$.
We used this fact for the second equality above;
 see also \cite{kl07} for a different argument leading to the same conclusion.

From (\ref{me4}) and (\ref{me5}), it then follows that
$$\begin{array}{ll} \tP^y[d(y,\xi_{T_1})\ge r] & = \dis\sum_{k\ge 1} \tP^y[d(y,\xi_{S_k})\ge r;\,S_k=T_1]\\
& \\
& = \dis\sum_{k\ge 1} {\widetilde \EE}^y[\tP^y[S_k=T_1 \vert\, \xi]\,;\, d(y,\xi_{S_k})\ge r]\\
& \\
& \lesssim \dis\sum_{k\ge 1} (1-(1/C^2))^k \tP^y[d(y,\xi_{S_k})\ge r]\\
& \\
& \lesssim e^{-\la r} \dis\sum_{k\ge 1} (1-(1/C^2))^ke^{\la c_4 k}\,.\end{array}$$
Thus, there is some $\la_0>0$ so that if we choose $\la\in(0,\la_0]$ then this last series is convergent and we find
$$ \tP^y[d(y,\xi_{T_1})\ge r] \lesssim e^{-\la r}\,.$$

Consequently, noting that $d(Z_1(w),\xi_{T_1})\le \ep$ and choosing $\la=\la_0$,
$$\EE\left[e^{(\la_0/2)d(y,Z_1(w))}\right] \lesssim 1 +\int_{u > 0} e^u\tP[d(y,\xi_{T_1}) > 2u/\la_0]du \lesssim 1+ \int_{u>0} e^{-u}du<\infty\,.$$
\endp

\subsection{Examples}\label{exples}
Let us fix $n\ge 2$ and consider the hyperbolic space $\HH^n$ of constant sectional curvature $-1$.
The explicit form of the Green function on this space shows easily that, given $w,x,y,z\in \HH^n$
which are at distance $c>0$ apart from one another, one has
\REFEQN{anchyp} \Te(x,y)\gtrsim \min\{\Te(x,z),\Te(z,y)\}\ENDEQN
 where $\Te$ is Na\"{\i}m's kernel, and
the implicit constant depends only on $c$.
Let $N$ be a finite volume hyperbolic manifold with deck transformation group $\G$ acting on $\HH^n$.
The estimate (\ref{anchyp}) shows that the Green
metric $d_G$ on $\G$ associated with the discretised Brownian motion on $\HH^n$ is hyperbolic.
Moreover, the estimate (ED) holds as well, so that the Green metric $d_G$
is quasi-isometric to the restriction of the hyperbolic metric to the orbit $\G(o)$ of a base point $o\in\HH^n$.
Since $N$ has finite volume, the limit
set of $\G$ is the whole sphere at infinity, and it coincides with the visual boundary of $(\G,d_G)$.
Therefore, \Th \ref{martin} implies that
the Martin boundary coincides with $\partial \HH^n$, homeomorphic to $\SS^{n-1}$.
We omit the details.

We apply this construction in  two special cases.

If we consider for $N$ a punctured $2$-torus with a complete hyperbolic metric of finite volume
(as in \cite{bl}), we obtain an example of a random walk on the free group for which the Green metric
is hyperbolic but its boundary $\SS^1$ does not coincide with the boundary of the group (which is a Cantor set).
Therefore, $d_G$ does not belong to the quasi-isometry class of the free group.

If we consider now for $N$ a complete hyperbolic $3$-manifold of finite volume with a rank $2$ cusp, then its
fundamental group is not hyperbolic since it contains a subgroup isomorphic to $\Z^2$, but the Green
metric is hyperbolic nonetheless.


\appendix

\section{Quasiruled hyperbolic spaces}\label{app}

For geodesic spaces, hyperbolicity admits many characterisations based on  geodesic triangles (cf. Prop.\,2.21 from \cite{gh}).
Most of them still hold when the space $X$ is just a length space (see eg. \cite{vag}).
For instance, a geodesic hyperbolic space satisfies Rips condition, namely, a constant $\de$ exists
such that any edge of a geodesic triangle is at distance at most $\de$ from the two other edges.

It is known that if $X$ and $Y$ are two quasi-isometric geodesic spaces, then
$X$ is hyperbolic if and only if $Y$ is (\Th 5.12 in \cite{gh}). This statement is known to be false in general if
we do not assume both spaces to be geodesic (Example 5.12 from \cite{gh}, and Proposition \ref{nonhyp} below).

Since quasi-isometries do not preserve small-scales of metric spaces, in particular
geodesics, it is therefore important to find other coarse characterisations of hyperbolicity.
Such a characterisation is the purpose of this appendix. We propose a setting which enables
us to go through the whole theory of \qc measures as if the underlying space was geodesic.

\noindent{\bf Definition.} A {\it quasigeodesic curve (resp. ray, segment)} is the image of
$\R$ (resp. $\R_+$, a compact interval of $\R$) by a quasi-isometric embedding. A space is said
to be {\it quasigeodesic} if there are constants $\la$, $c$ such that any pair of points
can be connected by a $(\la,c)$-quasigeodesic.

The image of a geodesic space by a quasi-isometry is thus quasigeodesic. But as it was mentioned earlier,
hyperbolicity need not be preserved.

\noindent{\bf Definition.} A {\it $\tau$-quasiruler} is a quasigeodesic
$g:\R\to X$ (resp. quasisegment $g:I\to X$,  quasiray $g:\R_+\to X$) such that, for
any $s<t<u$, $$(g(s)|g(u))_{g(t)}\le \tau.$$

Let $X$ be a metric space. Let $\la \ge 1$ and $\tau,c>0$ be constants.
A {\it quasiruling structure} $\GGG$  is
 a set of $\tau$-quasiruled $(\la,c)$-quasigeodesics such any
pair of points of $X$ can be joined by an element of $\GGG$.

A metric space will be {\it quasiruled} if constants $(\la,c,\tau)$ exist so that the space is
$(\la,c)$-quasigeodesic and if every $(\la,c)$-quasigeodesic is
a $\tau$-quasiruler i.e., the set of quasigeodesics defines a quasiruling structure.
The data of a quasiruled space are thus the constants $(\la,c)$ for the quasigeodesics
and the constant $\tau$ given by the quasiruler property of the $(\la,c)$-quasigeodesics.

A quasi-isometric embedding $f:X\to Y$ between a geodesic metric space $X$ into a metric space
$Y$ is {\it $\tau$-ruling} if the image of any geodesic segment is a $\tau$-quasiruler.
 Then
the images of geodesics of $X$ define a quasiruling structure $\GGG$ of $Y$.
In this situation, we will say that
$\GGG$ is induced by $X$.

\REFTHM{main} Let $X$ be a geodesic hyperbolic metric space, and $\vp:X\to Y$ a quasi-isometry, where $Y$ is
a metric space. The following statements are equivalent:
\ben
\item[(i)] $Y$ is hyperbolic;
\item[(ii)] $Y$ is quasiruled;
\item[(iii)]  $\vp$ is ruling.
\een

 Moreover if $Y$ is a hyperbolic
quasiruled space, then
 $Y$ is isometric to a quasiconvex subset of a geodesic hyperbolic metric space $Z$.

 Furthermore, if $\Gamma$ acts geometrically on $Y$, then $\Gamma$ is a quasiconvex group acting on $Z$.
\ENDTHM

\Th \ref{hypqrule} is a consequence from \Th \ref{main}.

We refer to \cite{gh} for any undefined notion used in the sequel.

\subsection{Straightening of configurations}

Let $I=[a,b]\subset \R$ be a closed connected subset.
We assume throughout
this section that constants $(\la,c,\tau)$ are fixed.

\REFLEM{strl1} Let $g:I\to X$ be a quasiruler. There is a $(1,c_1)$-quasi-isometry
$$f:g(I)\to [0, |g(b)-g(a)|]\,,$$ for some $c_1$ which depends only on the data
($\la$, $c$ and $\tau$).\ENDLEM

\Proof For any $x\in g(I)$, let $f(x)=\min\{|x-g(a)|,|g(b)-g(a)|\}$. Thus
\REFEQN{eq} ||x-g(a)|-f(x)|\le 2\tau. \ENDEQN

Let $x,y\in g(I)$ with $x=g(s)$ and $y=g(t)$, and let us assume that $s<t$.
\bit
\item  We apply (\ref{eq}) repeatedly.
On the one hand, $$|f(x)-f(y)|\le ||x-g(a)|-|y-g(a)||+4\tau\le |x-y|+4\tau\,.$$
On the other hand, since $s<t$, it follows that $$|x-g(a)|+|x-y|\le |y-g(a)|+2\tau$$ so that
$$|f(x)-f(y)|\ge  |x-y|-8\tau.$$
Hence $f$ is a $(1,8\tau)$-quasi-isometric embedding.

Note that the constants above are not sharp (a case by case treatment would divide most
of them by $2$).

\medskip

\item If $|a-b|\le 2$, then $|f(g(a))-f(g(b))|=|g(a)-g(b)|\le 2\la +c$ and $f$ is cobounded.

Otherwise, $|a-b|> 2$. Let $s_j=a+j$ for $j\in\N\cap [0,|b-a|]$. It follows that
$$|f(g(s_j))-f(g(s_{j+1}))|\le \la|s_j-s_{j+1}|+c+4\tau\le \la+c+4\tau.$$
The set $\{f(g(s_j))\}_j$ is a chain in $[0,|g(b)-g(a)|]$ which joins $0$ to $|f(g(a))-f(g(b))|=|g(b)-g(a)|$;
since two consecutive points of $\{f(g(s_j))\}_j$ are at most
$\la+c+4\tau$ apart, it follows that its  $(\la+c+4\tau)$-\nbhd covers  $[0,|g(a)-g(b)|]$,
hence $f$ is a quasi-isometry.\endp
\eit

\noindent{\bf Remark.} If $f_a$ denotes the map as above and $f_b:g(I)\to [0, |g(b)-g(a)|]$
the map such that $f_b(g(b))=0$, then
$|f_a(x)+f_b(x)- |g(a)-g(b)||\le 2\tau$ holds.

\noindent{\bf Definition.}
Given three points $\{x,y,z\}$, there is a tripod $T$ and an isometric embedding $f:\{x,y,z\}\to T$ such that
the images are the endpoints of $T$. We let $\bar c$ denote the center of $T$.

A quasitriangle $\De$ is given by three points $x,y,z$ together with three quasirulers
joining them. We will denote the edges by $[x,y]$, $[x,z]$ and $[y,z]$.
Such a quasitriangle is $\de$-thin if any segment is in the $\de$-\nbhd of the two others.

\REFLEM{strl2} Let $\De$ be a $\de$-thin quasitriangle with vertices $\{x,y,z\}$. There is a $(1,c_2)$-quasi-isometry
$$f_{\De}:\De\to T\,,$$ where $T$ is the tripod associated with $\{x,y,z\}$ and $c_2$ 
depends only on the data ($\delta$, $\la$, $c$, $\tau$).\ENDLEM

\Proof Let us define $f_{\De}$ using Lemma \ref{strl1} on each
edge. This map is clearly cobounded.

Let $u,v\in\De$. Since $\De$ is thin, one may find two points $u',v'\in\De$ on the same edge such
that $|u-u'|\le \de$ and $|v-v'|\le \de$, so that $$||u-v|-|u'-v'||\le 2\de.$$

If $u$ and $u'$ belong to the same edge, then
$$
|f_{\Delta} (u) -f_{\Delta} (u')|\leq |u-u'|+c_1\leq \delta +c_1 \, .
$$
Otherwise, let $x$ be the common vertex of the edges containing $u$ and $u'$,
then it follows from (\ref{eq})
that $$|f_{x}(u)-f_{x}(u')|\le |u-u'| +4 \tau\le \de +4 \tau$$ and similarly for $v$ and $v'$.
Thus $$|f_{\De}(u)-f_{\De}(u')|,|f_{\De}(v)-f_{\De}(v')| \le c',$$ where $c'$ depends only on the data.

It follows that $$||f_{\De}(u)-f_{\De}(v)|-|f_{\De}(u')-f_{\De}(v')||\le 2c'.$$

But since $u'$ and $v'$ belong to the same edge, Lemma \ref{strl1}
implies that $$||f_{\De}(u')-f_{\De}(v')|-|u'-v'||\le c_1\,,$$
so
$$||f_{\De}(u)-f_{\De}(v)|-|u'-v'||\le 2c'+c_1 $$ and finally
$$||f_{\De}(u)-f_{\De}(v)|-|u-v||\le (2c'+c_1+2\de).$$\endp

\bigskip

In the situation of Lemma \ref{strl2} we have $$|(f_{\De}(x)|f_{\De}(y))_{f_{\De}(z)}-(x|y)_{z}|\le C\,,$$
\noindent for some universal constant $C>0$;
thus,  we may find points
$c_x\in [y,z]$, $c_y\in [x,z]$ and $c_z\in [y,x]$ such that
$$|f_{\De}(c_x)-\bar c|,|f_{\De}(c_y)-\bar c|,|f_{\De}(c_z)-\bar c|\le c_3,$$
and $$\diam\{c_x,c_y,c_z\}\le c_3,$$ where
$c_3$ depends only on the data.

\REFPROP{strp} Let $X$ be a metric space endowed with a quasiruling structure $\GGG$
such that all quasitriangles are
$\de$-thin. Then $X$ is hyperbolic quantitatively: the constant of hyperbolicity only depends
on $(\de,\la,c,\tau)$.\ENDPROP

\Proof Let us fix $w,x,y,z\in X$. Let us consider the following triangles\,:
$A=\{w,x,z\}$ and $B=\{w,x,y\}$. 
Let us denote by $T_{A}$, $T_{B}$ and $\bar c_A$, $\bar c_B$
the associated tripod and center respectively, and let us define $Q=T_A\cup T_B$ where both copies  $f_A([w,x])$
and  $f_B([w,x])$ of
$[w,x]$ have been
identified. This metric space $Q$ is topologically an ``$\times$'', and so is of course $0$-hyperbolic.

Let us define $f:A\cup B\to Q$ by sending $A$ under $f_A$ and $B$ under $f_B$.

The restriction of  $f$ to $A$ and to $B$ is a $(1,c_2)$-quasi-isometry by Lemma \ref{strl2}.

It follows that
$$|f(y)-f(z)|=|f(y)-\bar c_B|+|\bar c_B-\bar c_A|+|\bar c_A-f(z)|.$$
One may find $c_A,c_B\in [w,x]$ such that $|f(c_A)-\bar c_A|\le c_3$ and $|f(c_B)-\bar c_B|\le c_3$.
Lemma \ref{strl2} implies that $|f(y)-f(c_B)|=|y-c_B|$ and $|f(c_A)-f(z)|=|c_A-z|$ up to an additive constant.
Therefore, $|f(y)-\bar c_B|=|y-c_B|$ and $|\bar c_A-f(z)|=|c_A-z|$ up to an additive constant too.
By Lemma \ref{strl1}, $|\bar c_B-\bar c_A|=|c_B-c_A|$ up to an additive constant, whence the existence
of some constant
$c_4>0$ such that
$$|f(y)-f(z)|\ge |y- c_B|+| c_B- c_A|+|c_A-z|-c_4\ge |y-z|-c_4.$$

Hence  $(f(y)|f(z))_{f(w)}\le (y|z)_{w}+c_4$. It follows from the hyperbolicity of $Q$ that
$$(y|z)_{w}  \ge \min\{(f(x)|f(z))_{f(w)}, (f(y)|f(x))_{f(w)}\} -c_4$$
and since the restrictions of $f$ to $A$ and $B$ are $(1,c_2)$-quasi-isometries,
$$ \min\{(f(x)|f(z))_{f(w)}, (f(y)|f(x))_{f(w)}\} -c_4\ge \min\{(x|z)_{w}, (y|x)_{w}\} -c_5$$
for some constant $c_5$. We have just established that for any $w,x,y,z$,
$$(y|z)_{w}  \ge\min\{(x|z)_{w}, (y|x)_{w}\} -c_5.$$\vspace*{-0.5cm}\endp

\subsection{Embeddings of hyperbolic spaces}

We recall a \th of M.\,Bonk and O.\,Schramm (\Th 4.1 in \cite{bsc})\,:

\REFTHM{bksch} Any $\de$-hyperbolic space $X$ can be isometrically embedded into a complete geodesic
$\de$-hyperbolic space $Y$.\ENDTHM

We will show that if $\Gamma$ acts isometrically on $X$, then so is the case on $Y$. To prove this we need to review
the construction of the set $Y$.

The first lemma, which we recall, is the basic step in the construction.

\REFLEM{bscl1} Let $X$ be $\de$-hyperbolic metric space, and let $a\ne b$ be in $X$. If, for every $x$,
$(|a-b|/2,|a-b|/2)\ne (|a-x|,|b-x|)$, then there is a $\de$-hyperbolic space $X[a,b]=X\cup\{m\}$ such that
$(|a-b|/2,|a-b|/2)= (|a-m|,|b-m|)$. Furthermore, for any $x\in X$,
$$|x-m|= \frac{|a-b|}{2} +\sup_{w\in X}(|x-w|-\max\{|a-w|,|b-w|\})\,.$$\ENDLEM

We call $m$ the middle point of $\{a,b\}$.

\REFLEM{bscl2} A $\de$-hyperbolic metric space $X$ embeds isometrically into a $\de$-hyperbolic space $X^*$ such that, for any $(a,b)\in X$,
there exists a middle point $m=m(a,b)\in X^*$. \ENDLEM

\Proof They apply a transfinite induction\,: let $\phi:\om \to X\times X$ be an ordinal of $X\times X$. Define inductively $X(\al)$ as follows.
Set $X(0)=X$. If $\al=\be+1\le \om+1$, then define $X(\al)=X(\be)[\phi(\al)]$. Clearly, $X(\al)$ is $\de$-hyperbolic.
If $\al$ is a limit ordinal, set $$X(\al)=\left(\cup_{\be <\al} X(\be)\right)[\phi(\al)]\,.$$
Here too, $X(\al)$ is $\de$-hyperbolic since $\de$-hyperbolicity is preserved under increasing unions.
The space $X^*=X(\om +1)$ fulfills the requirements.\endp

For $\al\le\om +1$, let us define $m(\al)=m(\phi(\al))$ the middle of $\phi(\al)=(a(\al),b(\al))$, and let $D(\al)=|a-b|$.
If $x^*\in X$, set $\al(x^*)=0$\,;
otherwise, let $P(x^*)$ be the set of ordinals $\al$ such that  $x^*\in X(\al)$. Let us define $\al(x^*)$ as the minimum of $P(x^*)$; it
follows that $x^*=m(\al)$. We let $D(x^*)=D(\al)$. We also write $\phi(\al)=(a(x^*),b(x^*))$.

\REFLEM{bscl3} Let $\al<\be$,  then
$$|m(\al)-m(\be)|=\frac{D(\be)}{2}+\sup_{w\in X(\al)}\{|w-m(\al)|-\max\{|w-a(\be)|,|w-b(\be)|\}\}\,.$$\ENDLEM

\Proof Let
\begin{eqnarray*}
&&Z=\\
&&\left\{\g\in\om,\ |m(\al)-m(\g)|=\frac{D(\g)}{2}
+  \sup_{w\in X(\al)}\{|w-m(\al)|-\max\{|w-a(\g)|,|w-b(\g)|\}\}\right\}\,.
\end{eqnarray*}
The set $Z$ contains $\{\g\le \al +1\}$ by definition. Let us assume that $\{\g<\be\}\subset Z$ for some $\be>\al$. Pick $\g\in Z$,
so that $\al<\g<\be$. Given $\ep >0$, there is some $w\in X(\al)$ so that
$$|m(\al)-m(\g)|\le \frac{D(\g)}{2}+ |w-m(\al)|-\max\{|w-a(\g)|,|w-b(\g)|\}+\ep\,.$$
Since $w\in X(\al)$ is fixed, $$|m(\g)-a(\be)|\ge \frac{D(\g)}{2} + |w-a(\be)|-\max\{|w-a(\g)|,|w-b(\g)|\}\,.$$
A similar statement holds for $b(\be)$ instead of $a(\be)$.
Therefore
\begin{eqnarray*}&&\max\{|m(\g)-a(\be)|,|m(\g)-b(\be)|\}
\ge\\ &&\frac{D(\g)}{2} + \max\{|w-a(\be)|,|w-b(\be)|\}-\max\{|w-a(\g)|,|w-b(\g)|\}\,,
\end{eqnarray*}
and
\begin{eqnarray*}
&&|m(\al)-m(\g)|-\max\{|m(\g)-a(\be)|,|m(\g)-b(\be)|\}
\le\\ && |m(\al)-w|- \max\{|w-a(\be)|,|w-b(\be)|\} +\ep \,.
\end{eqnarray*}
It follows that, for each $\al <\g<\be $ and $\ep>0$, there is some $w\in X(\al)$ such that the supremum in the definition of $|m(\al)-m(\be)|$
is attained within $X(\al)$ up to $\ep$.
Hence $\be\in Z$, so $Z=X^*$ by induction. \endp

\REFLEM{bscl4} Let $0<\al<\be$. Then $|m(\al)-m(\be)|$ can be computed as
$$\frac{D(\al)}{2}+\frac{D(\be)}{2}+\sup_{w,w'\in X}\{|w-w'|-(\max\{|w-a(\al)|,|w-b(\al)|\}+\max\{|w'-a(\be)|,|w'-b(\be)|\})\}\,.$$\ENDLEM

\Proof We endow $\om\times\om$ with the lexicographical order, and we consider $\om'=\{(\al,\be), \al<\be\}$.
We assume by transfinite induction that the lemma is true for any $(\al,\be)<(\hal,\hbe)$. By Lemma \ref{bscl3}, given $\ep>0$, there is some $\hw\in X(\hal)$
such that $$|m(\hal)-m(\hbe)|\le\frac{D(\hbe)}{2} +|\hw-m(\hal)|- \max\{|\hw-a(\hbe)|,|\hw-b(\hbe)|\}+\ep\,.$$
It follows from the induction assumption that there are points $w',w\in X$ such that
\begin{eqnarray*}
|m(\hal)-\hw|-\ep
& \le & \frac{D(\hal)}{2}+
\frac{D(\hw)}{2}+|w-w'|  - (\max\{|w-a(\hal)|,|w-b(\hal)|\}\\ && +\max\{|w'-a(\hw)|,|w'-b(\hw)|\})\,.
\end{eqnarray*}

But
\begin{eqnarray*}
\max\{|\hw-a(\hbe)|,|\hw-b(\hbe)|\}
&\ge &\frac{D(\hw)}{2}  +\max\{|w'-a(\hbe)|,|w'-b(\hbe)|\}\\ && -\max\{|w'-a(\hw)|,|w'-b(\hw)|\}\,,
\end{eqnarray*}
so
\begin{eqnarray*}
|m(\hal)-m(\hbe)|-2\ep
&\le &\frac{D(\hal)}{2}+\frac{D(\hbe)}{2} +
|w-w'|-(\max\{|w-a(\hal)|,|w-b(\hal)|\}\\ &&+\max\{|w'-a(\hbe)|,|w'-b(\hbe)|\})\,.
\end{eqnarray*}
This establishes the lemma.\endp

\REFCOR{bscc1} If $\Gamma$ acts on $X$ by isometry, then it acts also on $X^*$ by isometry. \ENDCOR

\Proof If $x^*\in X^*\setminus X$ and $g\in \Gamma$, we let $g(x^*)=m(g(a(x^*)),g(b(x^*))))$. The fact that $g:X^*\to X^*$ acts
by isometry follows from Lemma \ref{bscl4} since the distance between two points relies only on points inside $X$.\endp

The construction now goes as follows. Define $X_0= X$, and $X_{n+1}= X_n^*$, for $n\ge 0$. The space $X'=\cup_{n\in\N} X_n$ is
a metric $\de$-hyperbolic space such that any pair of points admits a midpoint in $X'$. Note that if $\Gamma$ acts on $X$ by
isometry, then it also acts by isometry on $X'$.

To obtain a complete geodesic space, M.\,Bonk and O.\,Schramm use again a transfinite induction. Let $\om_0$ be the first
uncountable ordinal. They define  a metric space $Z(\al)$ for each ordinal $\al<\om_0$ such that $Z(\al)\supset Z(\be)$ if $\al>\be$.
We set $Z(0)$ as the completion of $X'$. More generally, if $\al=\be+1$, define $Z(\al)$ as the completion of $Z(\be)'$.
For limit ordinals $\al$, we define $Z(\al)$ as the completion of $\cup_{\be<\al}Z(\be)'$. It follows that for each $\al<\om_0$,
the metric space $Z(\al)$ is complete, $\de$-hyperbolic, and admits an isometric action of $\Gamma$ if $X$ did.

The construction is completed by letting $Y=\cup_{\al<\om_0} Z(\al)$. As above, an action of a group $\Gamma$ by isometry on $X$ extends canonically
as an action by isometry on $Y$.

\subsection{Quasiruled spaces and hyperbolicity}

We prove \Th \ref{main}  in four steps.

\subsubsection{} Let us assume that $Y$ is a quasigeodesic $\de$-hyperbolic space. It follows from \Th \ref{bksch}
that there are a $\de$-hyperbolic geodesic metric space $\hat{Y}$ and an isometric embedding
$\iota:Y\to \hat{Y}$. Thus, for any quasigeodesic segment $g:[a,b]\to Y$, $\iota(g)$ shadows a genuine geodesic
$\hat{g}=[\iota(g(a)),\iota(g(b))]$
from $\hat{Y}$ at distance $H=H(\la,c,\de)$. In other words, for any $t\in [a,b]$, there is a point
$\hat{y}_t\in\hat{g}$ such that $|\iota(g(t))-\hat{y}_t|\le H$. It follows that
$$(g(a)|g(b))_{g(t)}\le (\iota(g(a))|\iota(g(b))_{\hat{y}_t}+H= H$$ since $\iota(g(a))$, $\iota(g(b))$
and $\hat{y}_t$ belong to a geodesic segment.

Therefore, $Y$ is quasiruled.

\subsubsection{} If $Y$ is quasiruled, then $\vp$ is ruling since the image
under $\vp$ is a quasigeodesic, hence a quasiruler by definition.

\subsubsection{} Let us now assume that $X$ is a geodesic hyperbolic space and $\vp:X\to Y$ is
a quasi-isometry into a metric space $Y$. It follows that $Y$ is quasigeodesic and that the edge of
the image
of any geodesic triangle is at a bounded distance from the two other edges i.e., triangles are $\de$-thin.
If $\vp$ is ruling, then Proposition \ref{strp} applies, and proves that $Y$ is hyperbolic.

\subsubsection{} The statement concerning group actions follows from above and the previous section.

\subsection{Non-hyperbolic invariant metric on a hyperbolic group} \label{nonhypmetric}

In \cite{gh}, the authors provide an example of a non-hyperbolic metric space quasi-isometric to $\R$. One could wonder if, in the case of groups, the
invariance of hyperbolicity holds for quasi-isometric and invariant metrics. In this section, we disprove this statement.

\REFPROP{nonhyp} For any hyperbolic group, a left-invariant metric
quasi-isometric to a word metric exists which is not hyperbolic. \ENDPROP

We are grateful to C.\,Pittet and I.\,Mineyev for having pointed out to
us the metric $d$ in the following proof as a possible candidate.

\Proof Let $\Gamma$ be a  hyperbolic group and let  $|.|$ denote a word metric. We define
the metric  $$d(x,y)=|x-y|+ \log (1+|x-y|)\,.$$

Clearly,  $|x-y|\le d(x,y)\le 2|x-y|$ holds and $d$ is left-invariant by $\Gamma$.

Let us prove that $(\Gamma,d)$ is not quasiruled, hence not hyperbolic by \Th \ref{main}.

Let $g$ be a geodesic for $|.|$ which we identify with $\Z$. 
Since $(\Gamma,d)$ is bi-Lipschitz to $(\Gamma,|\cdot|)$,
it is a $(2,0)$-quasigeodesic for $d$. But
$$d(0,n)+d(n,2n)-d(0,2n)= \log (1+n)^2/(1+2n)$$ asymptotically behaves as  $\log n$. Therefore  $g$
is not quasiruled.\endp

\section{Approximate trees and shadows}\label{appsha}
Approximate trees is an important tool to understand hyperbolicity in geodesic spaces. Here, we adapt
their existence to the setting of hyperbolic quasiruled metric spaces following E.\,Ghys and P.\,de la Harpe
(\Th 2.12 in \cite{gh}).

\REFTHM{aptrees}
Let $(X,w)$ be a $\de$-hyperbolic metric space and let $k\ge 0$.
\begin{itemize}
\item[(i)]If $|X|\le 2^k+2$, then there is a finite metric pointed tree $T$ and a map $\phi:X\to T$ such that\,:\\
\indent $\to\ \forall x\in X$, $|\phi(x)-\phi(w)|=|x-w|$\,,\\
\indent $\to\ \forall x,y\in X$, $|x-y|-2k\de\le |\phi(x)-\phi(y)|\le |x-y|$.
\item[(ii)] If there are $\tau$-quasiruled rays $(X_i,w_i)_{1\le i\le n}$ with $n\le 2^k$ such that
$X=\cup X_i$, then  there is a pointed $\R$-tree $T$ and a map
$\phi:X\to T$ such that\\
\indent $\to\ \forall x\in X$, $|\phi(x)-\phi(w)|=|x-w|$\,,\\
\indent $\to\ \forall x,y\in X$, $|x-y|-2(k+1)\de-4c -2\tau\le |\phi(x)-\phi(y)|\le |x-y|$,
where $c=\max\{|w-w_i|\}$.
\end{itemize}\ENDTHM

We repeat the arguments in \cite{gh}. The proofs of the first two lemmata can be found in \cite{gh},
and the last one is the quasiruled version of \cite[Lem. 2.14]{gh}. In the three lemmata, $X$ is assumed to
be $\de$-hyperbolic. Furthermore, we will omit the subscript $w$ for the inner product and
write $(\cdot|\cdot)=(\cdot|\cdot)_w$.

\REFLEM{aa1} We define
\begin{itemize}
\item[$\to$]  $(x|y)'=\sup\min\{(x_{i-1}|x_i),\ 2\le i\le L\}$, where the supremum is taken over
all finite chains $x_1,\ldots,x_L$ with $x_1=x$
and $x_L=y$,
\item[$\to$] $|x-y|'=|x-w|+|y-w|-2(x|y)'$,
\item[$\to$] $x\sim y$ if $|x-y|'=0$.\end{itemize}
Then $\sim$ is an equivalence relation and  $|\cdot |'$ is a distance on $X/\sim$ which makes
it a $0$-hyperbolic space. Moreover,
for any $x\in X$, $|x-w|'=|x-w|$ holds, and for any $x,y\in X$, $|x-y|'\le|x-y|$.\ENDLEM

\REFLEM{aa2} If $|X|\le 2^k+2$ then for any chain $x_1,\ldots,x_L\in X$,
 $$(x_1|x_L)\ge \min_{2\le j\le L}\{(x_{j-1}|x_j)\}-k\de\,,$$ holds.
\ENDLEM

\REFLEM{aa3} Let $X=\cup_{i=1}^n X_i$ where $(X_i,w_i)$ are $\tau$-quasiruled rays. If $n\le 2^k$ then,
for any chain $x_1,\ldots,x_L\in X$,
$$(x_1|x_L)\ge\min_{2\le j\le L} \{(x_{j-1}|x_j)\}-(k+1)\de -2c-\tau\,.$$\ENDLEM

\Proof First, $(x|y)_w\le \min\{|x-w|,|y-w|\}$ holds for any $x,y\in X$, and if $x,y\in X_i$ then
$|(x|y)_{w_i}-\min\{|x-w_i|,|y-w_i|\}|\le \tau$, and
$|x-w_i|\ge |x-w|-|w-w_i|\ge|x-w|-c$. Similarly, $|y-w_i|\ge |y-w|-c$. Thus,
$(x|y)_{w_i}\ge \min\{|x-w|,|y-w|\} -c-\tau$ and
$$(x|y)_w\ge (x|y)_{w_i}-c\ge \min\{|x-w|,|y-w|\} -2c-\tau\ge\min\{ (x|x')_w, (y|y')_w\} -2c-\tau$$
for all $x',y'\in X$.

Let $x_1,\ldots,x_L\in X$ be a chain. We will write $X(x_j)$ to denote the quasiruled ray $X_i$ which
contains $x_j$. Either, for all $j\ge 2$, $x_j\not\in X(x_1)$, or there is
a maximal index $j>1$ such that
$x_j\in X(x_1)$. Hence, it follows from above that $(x_1|x_j)\ge \min_{2\le i\le j}\{(x_{i-1}|x_i)\}-2c-\tau$.
In this case, let us consider  $x_1,x_j, x_{j+1}, \ldots,x_L$.

We inductively extract a chain $(x_i')$ of length at most $2n\le 2^{k+1}$ which
contains $x_1$ and $x_L$ and such that at most two elements belong to a common $X_i$,
and in this case, they have successive indices. It follows from Lemma \ref{aa2} and from above
that
$$(x_1|x_L)\ge \min\{(x_{i-1}'|x_i')\}-(k+1)\de \ge \min\{(x_{i-1}|x_i)\}-(k+1)\de -2c -\tau\,.$$ \endp

{\noindent\sc Proof of \Th \ref{aptrees}.} The \th follows as soon as we have found a quasi-isometric
embedding $\phi:X\to T$ with $T$ $0$-hyperbolic.

Lemma \ref{aa1} implies that $X/\sim$ is $0$-hyperbolic  and that $\phi:X\to X/\sim$ satisfies $|\phi(x)-\phi(w)|'=|x-w|$ and
$|\phi(x)-\phi(y)|'\le|x-y|$.

For case (i), Lemme \ref{aa2} shows that $(x|y)\ge (x|y)'-k\de$ i.e.,
$$|\phi(x)-\phi(y)|'\ge|x-y| -2k\de.$$

For case (ii), Lemme \ref{aa3} shows that  $(x|y)\ge (x|y)'-(k+1)\de-2c-\tau$ i.e.,
$$|\phi(x)-\phi(y)|'\ge|x-y| -2(k+1)\de -4c-2\tau.$$\endp

\bigskip

{\noindent\bf Visual quasiruling structures.}
Let $(X,d,w)$ be a hyperbolic space endowed with a quasiruling structure
$\GGG$. We say that $\GGG$ is {\it visual} if any pair of points in $X\cup\partial X$ can be joined by
a $\tau$-quasiruled $(\la,c)$-quasigeodesic. If $X$ is a proper space, then any quasiruling structure
can be completed into a visual quasiruling structure. Also, if $Y$ is a hyperbolic geodesic proper metric
space and $\vp:Y\to X$ is ruling, then the induced quasiruling structure is also visual.
 This fact can in particular be applied when $Y$ is a locally finite Cayley graph of a non-elementary hyperbolic group
$\G$, $(X,d)\in\DDD(\G)$ and $\vp$ is the identity map. Thus one endows $(X,d)$ with a
visual quasiruling structure.

\bigskip

{\noindent\bf Shadows.} Let $(X,d,w)$ be a hyperbolic  quasiruled space endowed with a visual quasiruling structure
$\GGG$.
We already defined the shadow $\mho(x,R)$ in Section 2  as the set of points $a\in\partial X$ such
that $(a|x)_w\ge d(w,x)-R$.  An alternative definition is:
let  $\mho_{\GGG}(x,R)$ be the set of points $a\in\partial X$ such
that there is a quasiruler $[w,a)\in\GGG$ which intersects $$B(x,R) = \{y \in X\, :\ d(x,y)< R\}.$$
The following holds by applying  \Th \ref{aptrees}, since  $\GGG$ is visual.

\REFPROP{shadowVSballbis}Let $X$ be a hyperbolic space endowed with a visual quasiruling structure $\GGG$.
There exist positive constants $C,R_0$ such that for any $R>R_0$, $a\in \partial X$ and $x\in [w ,a)\in\GGG$,
$$
 \mho_{\GGG} (x,R-C)\subset \mho (x,R) \subset
\mho_{\GGG} (x,R+C)  \, .
$$

\ENDPROP

The whole  theory of \qc measures for hyperbolic groups acting on geodesic spaces in \cite{co} is
based on the existence of approximate trees. Therefore, the same proof as in \cite{co} leads to \Th \ref{reg2}
and Lemma \ref{shadow}. Since \qc measures are Ahlfors-regular, the lemma of the shadow also holds for
shadows defined by visual quasiruling structures.

 Note that,  in a hyperbolic space endowed with a visual quasiruling structure,
\Th \ref{aptrees} implies that
the definition of Busemann functions we gave in Section \ref{S2} is equivalent to  the classical one
given below:

{\noindent\bf Busemann functions.} Let us assume that $(X,w)$ is a pointed hyperbolic quasiruled space.
Let $a\in\partial X$, $x,y\in X$ and $h:\R_+\to X$ a quasiruled ray such that  $h(0)=y$ and $\lim_{\infty}h=a$.
We define $\be_a(x,h) \fdefeq \limsup (|x-h(t)|-|y-h(t)|)$ and
$$\be_a(x,y) \fdefeq \sup\{\be_a(x,h),\mbox{ with }h\mbox{ as above}\}\,.$$

One can actually retrieve the metric from the Busemann functions as the next Lemma shows.  

\REFLEM{aa4}
Let $(X,w)$ be a pointed hyperbolic quasiruled space with the following 
quasi-starlike property:  there exists $R_1$ such that any $x\in X$ is at distance at most $R_1$ from 
a quasiray $[w,a)$, $a\in\partial X$. 
Then there exists a constant $c_6$ such that 
$$ \vert \vert x-y\vert -\sup_{a\in\partial X} \beta_a(x,y)\vert\leq c_6\,,$$ 
for all $x,y\in X$.  

The constant $c_6$ depends only on the data ($\delta$, $\la$, $c$, $\tau$). 
\ENDLEM

Observe that the quasi-starlike property is satisfied as soon as there is a geometric group 
action on $X$. 

\Proof From the triangle inequality, we always have $\be_a(x,y)\leq \vert x-y\vert$. 
Now choose $x,y\in X$ and $a\in\partial X$ such that $y$ is at distance at most $R_1$ 
from  a quasiray $[w,a)$. If these four points were really sitting on a tree, we would have 
$\be_a(x,y)\geq \vert x-y\vert -R_1$. Using approximate trees as in Theorem \ref{aptrees}, 
we get the Lemma. \endp


\end{document}